\documentclass[11pt]{article}
\usepackage[english]{babel}
\usepackage{amssymb,amsmath}
\usepackage{subcaption,graphicx}
\usepackage{amsthm}
\usepackage{float}
\usepackage{afterpage}
\usepackage{color}
\usepackage[round,authoryear]{natbib}
\usepackage{enumitem}
\usepackage{setspace}
\usepackage[table]{xcolor}
\usepackage{multirow}
\usepackage{natbib}
\usepackage{microtype}
\usepackage{enumitem}
\usepackage[textsize=scriptsize]{todonotes}
\usepackage[hidelinks]{hyperref}
\usepackage[ruled,vlined,linesnumbered]{algorithm2e}
\usepackage{booktabs}
\usepackage{pdflscape}
\usepackage{soul}
\usepackage[]{placeins}
\definecolor{light-gray}{gray}{0.95}
\sethlcolor{light-gray}
\SetKw{Continue}{continue}
\SetAlCapSkip{0.5em}
\DontPrintSemicolon
\SetKwInput{KwData}{Input}
\SetKwInput{KwResult}{Output}
\SetKwFor{For}{for}{}{endfor}
  \usepackage{algorithmic}
  \usepackage{longtable}
\newcolumntype{H}{>{\setbox0=\hbox\bgroup}c<{\egroup}@{}}

  \bibpunct[, ]{(}{)}{,}{a}{}{,}%
  \def\newblock{\ }%

\theoremstyle{definition}
  \newtheorem{definition}{Definition}

  \newtheorem{property}{Property}

  \newcommand{\cN}{{\mathcal{N}}}
  \newcommand{\cO}{{\mathcal{O}}}
  \newcommand{\cS}{{\mathcal{S}}}
\newcommand{\Bk}{{B$^k$}} 
\newcommand{\BS}{{\textnormal{B{\scriptsize{\&}}S}}}
\newcommand{\SA}{\mathcal{S}^{\textnormal{A}}}
\newcommand{\SB}[1]{\mathcal{S}^{\textnormal{B}}_{#1}}

    {\begin{list}{}{%
    \settowidth{\labelwidth}{\hspace{0cm}\textbf{#1} :}%
    \setlength{\leftmargin}{\labelwidth}\addtolength{\leftmargin}{\labelsep}}}%
    {\end{list}}

  \makeatletter
  \def\input@path{{input/}{tex/}{}}
  \makeatother
  \graphicspath{{image/}{}}%

  \definecolor{Emerald}{rgb}{0.2, 0.6, 0.3}
  \let\oldnl\nl
  \newcommand{\nonl}{\renewcommand{\nl}{\let\nl\oldnl}}
  
  \SetCommentSty{algcommentfont}

\begin{document}

\begin{center}
  \begin{Huge}
    Heuristics for vehicle routing problems: \vspace*{0.25cm} \linebreak  Sequence or set optimization? 
  \end{Huge}

  \vspace*{01cm}

  \textbf{T\'ulio A. M. Toffolo$^{1,2}$, Thibaut Vidal$^{3}$, Tony Wauters$^{1}$} \\
  
  \vspace*{0.2cm}
  {\small
      $^1$ KU Leuven, Department of Computer Science, CODeS \& imec - Belgium \\
      $^2$ Federal University of Ouro Preto, Department of Computing - Brazil \\
      $^3$ Pontif\'{i}cia Universidade Cat\'{o}lica do Rio de Janeiro, Computer Science Department - Brazil \\
  }

  \vspace*{0.2cm}
\end{center}

\noindent
\textbf{Abstract.} We investigate a structural decomposition for the capacitated vehicle routing problem (CVRP) based on vehicle-to-customer ``assignment'' and visits ``sequencing'' decision variables. We show that an heuristic search focused on assignment decisions with a systematic optimal choice of sequences (using Concorde TSP solver) during each move evaluation is promising but requires a prohibitive computational effort. We therefore introduce an intermediate search space, based on the dynamic programming procedure of Balas \& Simonetti, which finds a good compromise between intensification and computational efficiency. A variety of speed-up techniques are proposed for a fast exploration: neighborhood reductions, dynamic move filters, memory structures, and concatenation techniques. Finally, a tunneling strategy is designed to \emph{reshape} the search space as the algorithm progresses.

The combination of these techniques within a classical local search, as well as in the unified hybrid genetic search (UHGS) leads to significant improvements of solution accuracy. New best solutions are found for surprisingly small instances with as few as 256 customers. These solutions had not been attained up to now with classic neighborhoods. Overall, this research permits to better evaluate the respective impact of sequence and assignment optimization, proposes new ways of combining the optimization of these two decision sets, and opens promising research perspectices for the CVRP and its variants.\\


\section{Introduction} \label{sec:introduction}

The Capacitated Vehicle Routing Problem (CVRP) is classically described as a combination of a Traveling Salesman Problem (TSP) with an additional capacity constraint which lends a Bin Packing (BP) substructure to the problem \citep{Toth2014}. It can be seen as a Set Packing (SP) problem in which the cost of each set corresponds to the distance of the associated optimal TSP tour \citep{Balinski1964}. These problem representations emphasize the two decision sets at play: customer-to-vehicle \textsc{Assignments}, and \textsc{\textsc{Sequencing}} choices for each route \citep{Vidal2012a}, a duality that has left long-standing impressions in the literature, from the early developments of route-first cluster-second \citep{Bodin1979,Beasley1983} and cluster-first route-second constructive methods \citep{Fisher1981}, all the way to the set-covering-based exact methods and matheuristics which are currently gaining popularity.

Examining the recent progress on metaheuristics for the CVRP, little has changed in recent years concerning intra-route neighborhood search: \textsc{Relocate}, \textsc{Swap} and \textsc{2-opt}  neighborhoods and their immediate generalizations are employed, and these neighborhoods alone are sufficient to guarantee that most solutions resulting from a local search contain TSP--optimal routes. This is generally because classical CVRP instances involve short routes with up to 15 or 20 visits. For such small problems, even simple neighborhood search methods for the TSP tend to produce optimal tours.

  Based on this observation, a larger effort dedicated to TSP tour optimization, as a stand-alone neighborhood, is unlikely to result in further improvements. For this reason, it is very uncommon to observe the use of larger \emph{intra-route} neighborhoods (e.g., \textsc{3-Opt} or beyond) in recent state-of-the-art metaheuristics.
  Nevertheless, does this mean that \textsc{Sequencing} optimization should be abandoned in favor of more extensive search concerning \textsc{Assignment} choices? Certainly not. Indeed, even if local minima exhibit optimal TSP tours, inter-route moves frequently lead to TSP-suboptimal tours which are rejected due to their higher cost, but would be accepted otherwise if the tours were optimized. Such solution improvements would then not arise from separate \textsc{Assignment} or \textsc{Sequencing} optimizations, but from a careful combination of both.

\begin{figure}[!ht]
    \centering
    \vspace*{0.2cm}
    \includegraphics[width=0.65\textwidth]{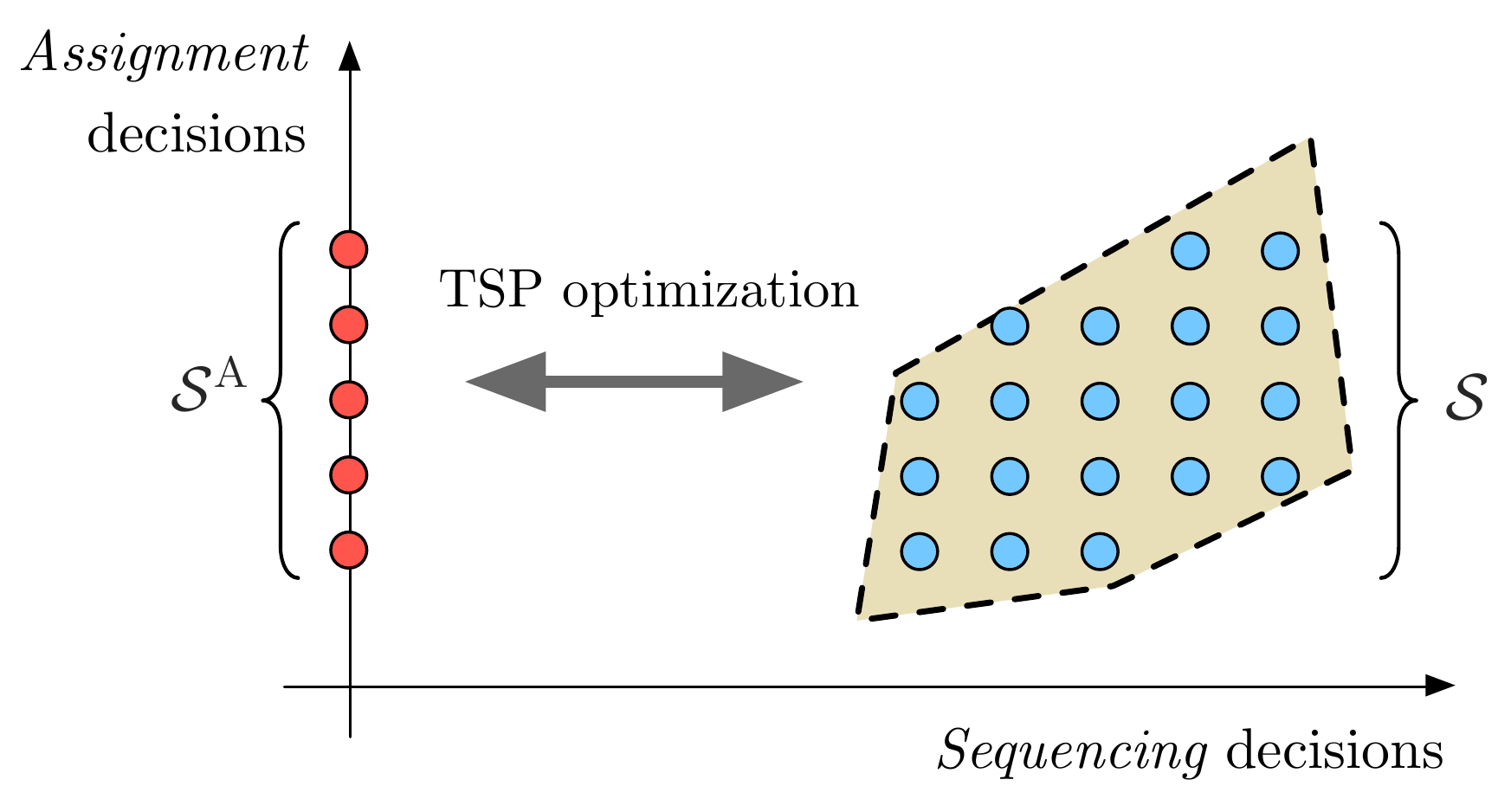}
    \caption{Two alternative search spaces for the CVRP}
    \label{fig:decision_projection}
  \end{figure}

Figure \ref{fig:decision_projection} schematically represents the solution set of the CVRP, whose decision variables are split into \textsc{Sequencing} decisions ($x$-axis) and \textsc{Assignment} decisions ($y$-axis). The $y$-axis also represents the solutions in terms of their \textsc{Assignment} decisions solely, ignoring \textsc{Sequencing} choices. These partial solutions can be viewed as a projection \citep{Geoffrion1970} of the original solutions $\cS$ on the space $\SA$ defined by a single decision subset (\textsc{Assignment}). Moreover, from a solution represented in terms of \textsc{Assignment} decisions, it is possible to find the best associated complete solution by solving each TSP associated with the routes.
 
  With this picture in mind, it is tempting to conduct a search in the space $\SA$ rather than~$\cS$. After all, the size of $\SA$ is exponentially smaller than that of $\cS$, the solutions of $\SA$ are in average of better quality, and the average size of a path in $\SA$ is smaller, such that fewer LS iterations are expected for convergence. However, the obvious drawback is that each move evaluation in $\SA$ requires solving one or several small TSPs to optimality, leading to a significant computational effort. Still, note that considerable progress has been made in the past~30 years with regard to the efficient solution of TSPs, and small problems with approximately 20 customers are solvable in a few milliseconds.
  Based on these observations, this article takes a fresh look at heuristic searches for the CVRP to answer two essential questions about the search space $\SA$:
  \begin{enumerate}[nosep]
    \item Is it practical and worthwhile to search in the space $\SA$ rather than $\cS$ ?
    \item If searching in $\SA$ requires an excessive effort, can we define a search space which maintains most of the key properties of $\SA$ but can be more efficiently explored~?
  \end{enumerate}

 As will be demonstrated in Section \ref{sec:experiments}, our experiments led us to answer the first question negatively: even with non-trivial memory and speedup techniques (hashtables and move filters) the computational overhead related to the exact solution of TSPs during each move evaluation, for a complete search in $\SA$, does not appear worth the gain in terms of solution quality.

By contrast, our answer to the second question is positive. Rather than requiring a complete exact solution of each TSP, the dynamic programming approach of \cite{Balas2001}, hereafter referred to as $\BS$, can be employed to perform a restricted route optimization \emph{during move evaluations}. Given a range parameter $k$ and an initial tour, the $\BS$ algorithm finds, in $\cO(k^2 2^{k-2}n)$ operations, the vertex sequence with minimum cost such that no vertex is displaced by more than $k$ positions. This allows us to define a search space $\SB{k}$ such that $\cS_0^{\textsc{B}} = \cS$ and $\lim\limits_{k \rightarrow \infty} \SB{k} = \SA$.
  Moreover, even for a fixed $k$, we propose \emph{tunneling} techniques that exploit the memory of past solutions to dynamically reshape the search space, in such a way that $\SB{k}$ converges towards $\SA$ as the search progresses.

To evaluate experimentally the potential of the new search spaces, we conduct experiments with a simple multi-start local search (MS-LS), and with the unified hybrid genetic search (UHGS) of \cite{Vidal2012,Vidal2012b}. The use of $\SB{k}$ for $k \in \{1,\dots,3\}$ appears to lead to solutions of higher quality on the new instances from \cite{Uchoa2017}. New best solutions were also found for surprisingly small instances with as few as 242 or 256 customers. These solutions had not been attained up to now with classic neighborhoods. Overall, this research allows to better evaluate the respective impact of \textsc{Sequencing} and \textsc{Assignment} optimization, proposing new ways to combine the optimization of these two decision sets, and leading to new state-of-the-art algorithms for the CVRP.

\section{Related Literature} \label{sec:literature}

  This section reviews some key milestones concerning the management of \textsc{Sequencing} and \textsc{Assignment} decisions in vehicle routing heuristics, as well as decision-set decompositions.

  \textsc{Sequencing} and \textsc{Assignment} decisions were first optimized separately in early constructive heuristics, giving rise to different families of methods. Route-first cluster-second algorithms \citep{Bodin1979,Beasley1983} first produce a giant TSP tour, before subsequently assigning consecutive visits into separate trips to produce a complete solution. In cluster-first route-second methods \citep{Fisher1981}, a clustering algorithm is employed to group customer visits into clusters, followed by TSP optimizations. Finally, \emph{petal algorithms} \citep{Foster1976,Renaud1996a} are based on an a-priori generation of candidate routes (petals), followed by the solution of a set packing or covering problem.

In the development of local search and metaheuristic algorithms which ensued in the 1990s and thereafter, \textsc{Assignment} and \textsc{Sequencing} optimizations began to be better integrated. 
Classical moves such as \textsc{Relocate}, \textsc{Swap}, \textsc{2-opt*} and their close variants allow to optimize both decision subsets. The associated neighborhood search methods form the basis of the vast majority of state-of-the-art algorithms. Petal algorithms have withstood the test of time, and high-quality routes are now extracted from local minima of metaheuristics instead of being enumerated in advance \citep[see, e.g., ][]{Muter2010,Subramanian2013}.

 The variety of vehicle routing problem variants has also triggered studies concerning problem decompositions.  \cite{Vidal2012a} established a review of the classical variants and their associated constraints, objectives, and decision sets, called \emph{attributes}. The attributes were classified in relation to their impact on \textsc{\textsc{Sequencing}}, \textsc{\textsc{Assignment}} decisions, and \textsc{Route evaluations} in heuristics, leading to a structural problem decomposition which serves as a basis for the Unified Hybrid Genetic Search (UHGS) algorithm of \cite{Vidal2012b} and allows to produce state-of-the-art results for dozens of VRP variants. Many problem attributes come jointly with new decision subsets, e.g., when optimizing vehicle routing with packing, timing or scheduling constraints \citep{Goel2012,Pollaris2015,Vidal2015b}, visit choices \citep{Vidal2014} or service-mode choices \citep{Vidal2014b,Vidal2015a}.

Decision-set decompositions are employed throughout many of the aforementioned papers to perform a search in the space of \textsc{Sequencing} and \textsc{Assignment} choices and optimally determine the remaining decision variables during each route and move evaluation.
 In this paper, the decision-set decomposition does not result from supplementary problem attributes, but is instead used to define exponential-size polynomially-searchable neighborhoods and transform the search space. Exponential-size neighborhoods have a long history in the combinatorial optimization literature \citep{Deineko2000,Ahuja2002,Bompadre2012}. Most of these neighborhoods are based on shortest path or matching subproblems, as well as specific graph and distance matrix structures with which some NP-hard problems become tractable (consider, as examples, Halin graphs or Monge matrices). As a rule of thumb, larger neighborhoods and faster search procedures are generally desirable. There are, however, theoretical limitations to the size of polynomially-searchable neighborhoods. \cite{Gutin2003} proved that, for the TSP, no neighborhood of cardinality at least $(n-k)!$ for a given constant $k$ can be searched unless NP $\subset$ P/poly.

  The neighborhood of \cite{Balas2001} is an exponential-size neighborhood for the TSP.
  Given an incumbent tour represented as a permutation $\sigma$ and a value $k$, it contains all permutations $\pi \circ \sigma$ such that $\pi$ fulfills $\pi(1) = 1$ and $\pi(i) \leq \pi(j)$ for all $i,j \in \{1,\dots,n\}$ such that $i+k \leq j$. In other words, if $i$ precedes $j$ by more than $k$ positions in $\sigma$, then $\pi(i)$ precedes $\pi(j)$. Setting $\pi(1) = 1$ allows to fix the origin location (e.g., depot). This neighborhood contains~$2^{\Theta(n)}$ solutions, and can be explored in $\cO(k^2 2^{k-2}n)$ using dynamic programming. This is a linear time complexity when $k$ is constant, and a polynomial complexity when $k = \cO(\log n)$. \cite{Balas2001} performed extensive experiments, and demonstrated that this dynamic programming procedure can be used as a stand-alone neighborhood to improve high quality local minima of the TSP and its immediate variants. Later, \cite{Irnich2008,Gschwind2016}, and \cite{Hintsch2017} employed this neighborhood to solve arc-routing problems with possible cluster constraints, and dial-a-ride problems.
One common characteristic of these studies is that they employed $\BS$ as a stand-alone neighborhood for route improvements. Only in one conference presentation \citep{Irnich2013}, the possibility of using the $\BS$ neighborhood \emph{in combination} with some classical CVRP moves has been highlighted, but the performance of such an approach remains largely unexplored.

We seek to go one step further.
Rather than applying this tour optimization procedure as a stand-alone optimization technique or in combination with a single classical neighborhood, we investigate its \emph{systematic use} in combination with every move of a classical CVRP local search. As discussed in the following, the methodological implications of such a redefinition of the search space are noteworthy.

\section{Proposed Methodology} 
\label{sec:methodology}

We will describe the methodology as a local search on indirect solution representations, using a decoder. This algorithm can be readily extended into a wide range of vehicle routing metaheuristics, e.g., tabu search, iterated local search, or hybrid genetic algorithm \citep{Gendreau2010e,Laporte2014a}.
There is no widely accepted term, in the current heuristic literature, for referring to the elements which represent such indirect solutions. The evolutionary literature usually refers to a \emph{genotype} to denote solution encodings (and \emph{phenotype} for the solutions themselves), whereas the local-search based metaheuristic literature refers to \emph{incomplete} or \emph{indirect} solutions (which are converted into complete solutions via a \emph{decoder} function). \emph{Incomplete} lets us think that the representation is necessarily a subset of a complete solution, and unnecessarily restricts the application scope. To circumvent this issue, we henceforth employ the term \emph{primitive solutions}. We first recall some basic definitions related to neighborhood search and indirect solution representations, then proceed with an analysis of alternative search spaces and the description of the proposed local search algorithm.

\begin{definition}[Primitive solutions and search space]
We consider a combinatorial optimization problem of the form $\min\limits_{x \in X} z(x)$, where $X$ is the solution space, and $z$ is an objective function to minimize. Let $Y$ be the set of primitive solutions, and let the decoder $f : Y \rightarrow X$ be an injective application that transforms any $y \in Y$ into a complete solution $x \in X$.
A neighborhood is defined as a mapping $\mathcal{N} : Y \rightarrow 2^Y$ that associates with each primitive solution $y$ a set of neighbors $\mathcal{N}(y) \subset Y$. The graph induced by $Y$ and $\mathcal{N}$ is referred to as the search space.
\end{definition}

\begin{definition}[TSP--optimal tour]
A tour $\sigma$ is TSP--optimal if there exists no other permutation of its visits $\pi \circ \sigma$ such that $\pi(1)=1$ with a shorter total distance.
\end{definition}

\begin{definition}[\Bk--optimal tour]
A tour $\sigma$ is \Bk--optimal if there exists no other permutation of its visits $\pi \circ \sigma$ with a shorter total distance such that $\pi(1)=1$ and $\pi(i) \leq \pi(j)$ for all $i,j \in \{1,\dots,n\}$ with $i+k \leq j$. The parameter $k$ is the range of the $\BS$ neighborhood.
\end{definition}

\subsection{A Choice of Search Space}
\label{sec:search_space}

In this section, we examine the search spaces associated with the set of all solutions ($\cS$), of those with TSP--optimal tours ($\SA$), and of those with \Bk--optimal tours ($\SB{k}$) and discuss their relative merits.\\

\noindent
\textbf{Search Space $\cS$.}
Classical local search methods for the CVRP do not distinguish between primitive and complete solutions. In the usual search space $\cS$, solution sets $X$ and $Y$ are equal and the decoder $f$ is the identity function. The neighborhood $\mathcal{N}$ is based on the definition of one or several classes of moves. A move $\phi$ is a local modification that can be applied to a primitive solution $y$ to generate a neighbor $\phi(y) \in \cN(y)$. 
For each search space considered in this paper, we will eventually refer to several classes of moves, but to a single neighborhood only, which corresponds to the union of all primitive solutions attainable from $y$ via one single move. Classical moves for the CVRP are based on relocations and exchanges of a bounded number of vertices, or replacements of a bounded number of edges. Most common neighborhoods have a quadratic cardinality ($|\mathcal{N}(y)| = \mathcal{O}(n^2)$ for all $y \in Y$). We refer to  \cite{Vidal2012a} for a comprehensive survey on classical local searches for the CVRP.

\begin{figure}[H]
      \centering
      \includegraphics[width=0.9\textwidth,page=1]{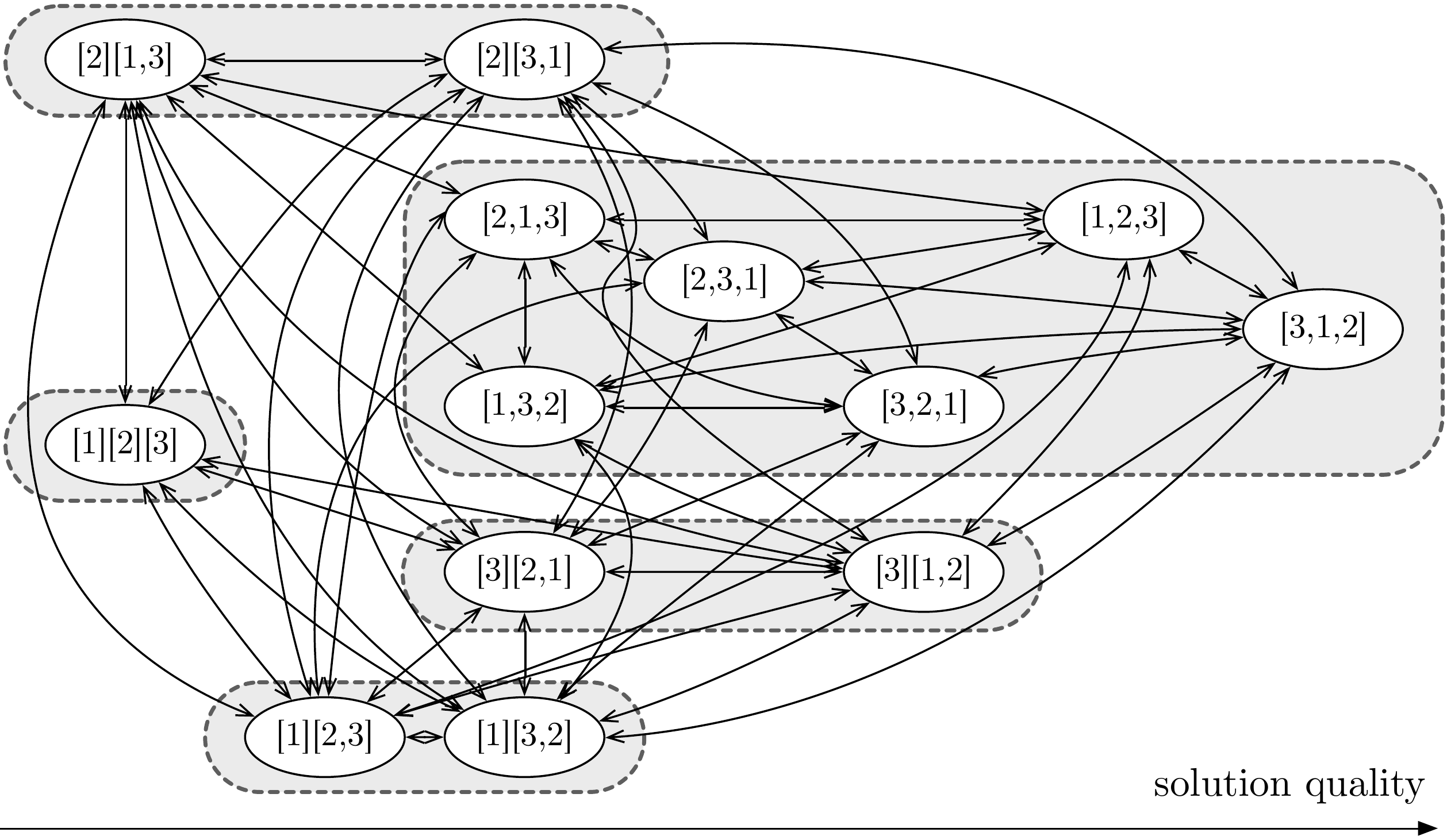}
      \vspace{-10pt}
      \caption{Search space $\cS$ for a small asymmetric CVRP instance} 
      \label{fig:neigh}
    \end{figure}

Figure \ref{fig:neigh} represents the search space $\cS$ associated with \textsc{Relocate} moves only, for a small asymmetric CVRP instance with three customers.
There are 13 possible solutions for this problem, each represented by a set of ordered customer visits. Solution `[1,2,3]', for example, employs one vehicle to visit customers 1, 2 and 3, while solution `[1][2][3]' employs three vehicles, one per customer. 
Each solution is represented by a node, positioned on the $x$-axis according to its quality (the more to the right, the better a solution is). The set of outgoing arcs of each solution points towards its neighbors.
Moreover, solutions with identical customer-route assignments are grouped within dashed areas. 
Note that, for this instance size, it is always possible to reach the optimum solution from any starting point in two successive moves.
In a local search that explores the neighborhood in random order and applies an improving move as soon as it is found, the worst case corresponds to five moves (when the initial solution is `[2][1,3]' or `[1][2][3]').\\

\noindent
\textbf{Search Space $\SA$.} As discussed earlier in this work, CVRP solutions can also be represented in terms of their \textsc{Assignment} decisions, excluding the \textsc{Sequencing} decisions in the representation and delegating the choices of the best visit sequences to the decoder. With such a paradigm, one can define a local search in the space $Y$ of primitive solutions, where each $y \in Y$ represents a partition of the customer set into subsets whose sums of demands do not exceed the vehicle capacity. The decoder $f$ is based on an exact TSP solver, responsible for generating the best visit sequence originating and finishing at the depot for each subset of customers. In this sense, the image $f[Y] \subset X$ contains exclusively solutions with TSP--optimal~tours.

\begin{figure}[H]
      \centering
      \includegraphics[width=0.9\textwidth,page=2]{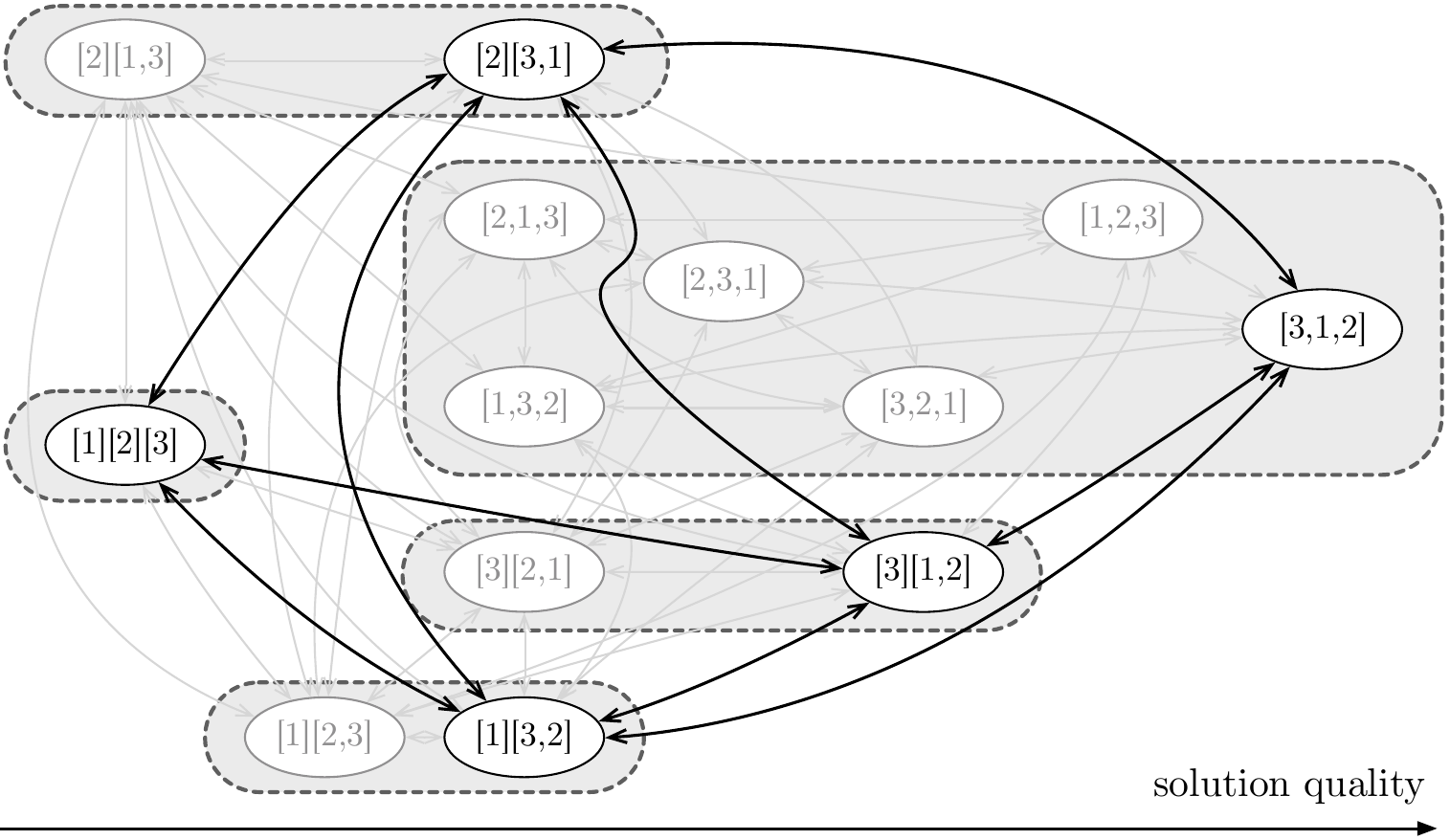}
      \vspace{-10pt}
      \caption{Search space $\SA$ for a small asymmetric CVRP instance} 
      \label{fig:neigh_concorde}
    \end{figure}

The neighborhood used to explore the search space $\SA$ can remain similar to classical CVRP neighborhoods, based on relocations or exchanges of customers between subsets, or involve other families of moves specialized for partition problems.
Figure \ref{fig:neigh_concorde} represents the resulting search space with simple \textsc{Relocate} moves.
Only TSP--optimal tours are explored and therefore the size of the search space reduces down to six primitive solutions.
The other solutions and their connections are represented in light gray.
Note, in our small example, that now at most three successive improving moves may be applied to attain the optimum from `[1][2][3]'.

Search space $\SA$ is smaller than $\cS$, and our computational experiments (Section~\ref{sec:experiments}) demonstrate that a search in this space indeed leads to solutions with higher quality. 
 However, each move evaluation in this space requires executing an algorithm with exponential worst-case time complexity, a TSP solver, in order to decode each primitive solution in the neighborhood for cost evaluation.
   Although research on the TSP has culminated in very efficient algorithms over the past thirty years, thousands (or millions) of small TSP instances should be solved during a local search in $\SA$, and thus the total computational effort dedicated towards decoding can grow prohibitively large. 
Moreover, bad behavior in a single case (e.g., due to an unusual long route with many customers or bad branching decisions) can be sufficient, without any other safeguard, to stall the entire algorithm.\\

\noindent
\textbf{Search Space $\SB{k}$.} To circumvent the aforementioned issues, we study an alternative search space in which the set of primitive solutions $Y$ is a subset of the complete solutions (with their \textsc{Assignment} and \textsc{Sequencing} decisions), but where the decoder $f$ is nontrivial, and consists of applying $\BS$ multiple times to each route with a fixed range ($k$ value) until the tours become \Bk--optimal. With these assumptions, the image $f[Y]$ contains exclusively complete solutions with \Bk--optimal tours.
           As such, the application of $\BS$ can be viewed as a post-optimization step \emph{during} classical CVRP move evaluations, opening the way for additional solution improvements. 
            A careful analysis of the resulting search space gives even more significance to this approach, due to three properties:

\begin{property}
\label{prop1}
From an initial solution containing a \Bk--optimal tour, a local search in the space~$\SB{k}$ explores only \Bk--optimal tours.
\end{property}

\begin{property}
\label{prop2}
For a fixed range $k$, each move evaluation and subsequent solution decoding is done in polynomial time as a function of $n$ and the number of applications of $\BS$.
\end{property}

\begin{property}
\label{prop3}
The search space $\SB{k}$ is such that $\SB{0} = \cS$ and $\SB{n-1} = \SA$, with $n$ being the number of customers.
\end{property}

These three properties are all fundamental for the methodology that follows.
Property~\ref{prop1} demonstrates how space $\SB{k}$ contains fewer solutions than $\cS$, and that the overall quality of these solutions tends to be higher (since non-\Bk-optimal tours are filtered out). Moreover, Property \ref{prop2} gives some computational time guarantees: even if the computational effort grows quickly with the range $k$, the effort of the decoder is guaranteed to remain stable for all $y \in Y$ when $k$ is constant, eliminating the possibility of a computational effort peak for specific TSP instances. Finally, Property \ref{prop3} demonstrates how $k$ balances the effort dedicated towards the optimization of the \textsc{Assignment} and \textsc{Sequencing} decision sets, and establishes $\SB{k}$ as an intermediate search space generalizing $\cS$ and $\SA$.

\begin{figure}[H]
      \centering
      \includegraphics[width=0.9\textwidth,page=3]{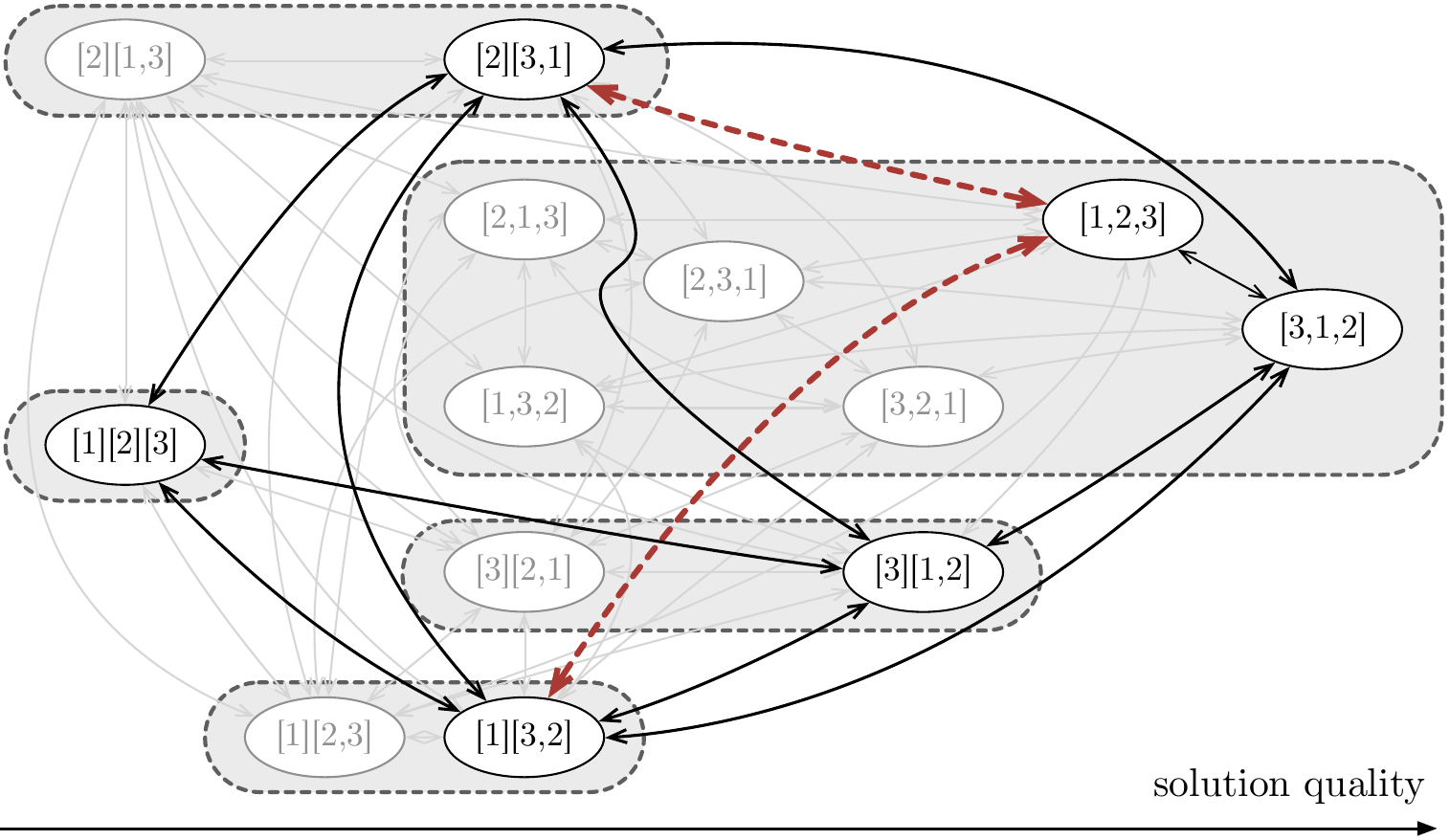}
      \caption{Search space $\SB{1}$ for a small asymmetric CVRP instance} 
      \label{fig:neigh_bs}
    \end{figure}

Figure \ref{fig:neigh_bs} illustrates the search space $\SB{k}$ for the same example as previous figures with $k=1$. It is an intermediate between the spaces depicted by Figures \ref{fig:neigh} and \ref{fig:neigh_concorde}, which correspond to $\SB{0} = \cS$ and $\cS^{\textsc{B}}_2 = \SA$, respectively. Note how the solution $[1,2,3]$ is now a neighbor of $[2][3,1]$ and $[1][3,2]$, as highlighted by the two dotted arcs.

\noindent
\textbf{Discussions and Choice.}
In light of these observations, we have conducted computational experiments on search space  $\SB{k}$, using the dynamic programming algorithm of \cite{Balas2001} to decode each solution, as well as on search space $\SA$ using the TSP solver \textsc{Concorde} \citep{Concorde2006}. Despite several speedup techniques (Section \ref{sec:lower_bound}), the search in space $\SA$ remained inefficient throughout our current experiments, especially for instances with a large number of customers per route. We therefore decided to focus on search space $\SB{k}$, and devised several speedup techniques to enable its efficient exploration.

\subsection{Efficient Local Search}
\label{sec:lower_bound}

To efficiently explore space $\SB{k}$, we developed a local search algorithm which exploits \emph{static neighborhood reductions}, \emph{dynamic move filters}, \emph{efficient memory structures} and \emph{concatenation techniques}. Most of these techniques seek to limit the search effort in~$\SB{k}$. This resulting method, displayed in Algorithm \ref{alg:overall}, can be easily integrated in state-of-the-art metaheuristics for vehicle routing problems.\\

     \begin{algorithm}[p]
       \caption{Efficient local search in the space $\SB{k}$}
      \label{alg:overall}
       \small
       \KwData{An initial complete solution $x^0$, an evaluation threshold $\psi$ and a granularity threshold $\Gamma$}

       $t \gets 0$

       \Repeat{\emph{$x^t$ is a local minimum}}
       {
\vspace*{0.15cm}
\tcp{Enumerating $\cO(\Gamma n)$ moves - candidate lists based on vertex proximity}
         \For{\emph{each move $\phi(x^t) \in \mathcal{N}(x^t)$ involving a vertex pair $(i,j)$, $j \in \Gamma(i)$}}
         {

\vspace*{0.15cm}
           \emph{The move $\phi$ modifies up to two routes of $x^t$. Let $z_\textsc{before}$ be the sum of the costs of these two routes, and let $(\sigma^1_1,\dots,\sigma^1_{b_1})$ and $(\sigma^2_1,\dots,\sigma^2_{b_2})$ be the new routes in $\phi(x^t)$.
}\;

	   \vspace*{0.3cm}
	  \tcp{First, filter infeasible moves with respect to capacity constraints~in~$\cO(1)$:}

	\vspace*{0.05cm}
           \If{$Q(\sigma^1_1 \oplus \dots \oplus \sigma^1_{b_1}) > Q$ \textbf{\emph{or}} $ Q(\sigma^2_1 \oplus \dots \oplus \sigma^2_{b_2}) > Q$ }
           {\vspace*{0.05cm} \Continue.}
          
          \vspace*{0.3cm}
          \tcp{Second, consider the cost of the classical CVRP move to filter non-promising solutions~in~$\cO(1)$:}

\vspace*{0.05cm}
         \If{$z(x^t) + C(\sigma^1_1 \oplus \dots \oplus \sigma^1_{b_1}) + C(\sigma^2_1 \oplus \dots \oplus \sigma^2_{b_2}) -  z_\textsc{before} >  (1 + \psi) \times z(x^t) $}
           {\vspace*{0.05cm} \Continue.}

           \vspace*{0.3cm}
             \tcp{Third, \emph{decode} the routes $\sigma^1$ and $\sigma^2$ to evaluate the move $\phi$ in  $\SB{k}$:}

 $z_\textsc{move} \gets 0$

	    \For{\emph{each route $\sigma^i$ with $i \in \{1,2\}$}}
{
			  \vspace*{0.15cm}
			\tcp{Compute hash key in~$\cO(1)$ and check memory in~$\cO(1)$:}

			$(\bar{\sigma}^i, \bar{z}_i) \gets \textsc{Lookup}(H(\sigma^i_1 \oplus \dots \oplus \sigma^i_{b_i}))$
			
			\vspace*{0.15cm}
			\tcp{If not in memory, repeatedly apply the dynamic programming algorithm of \cite{Balas2001} until the route becomes :}
			\If{$(\bar{\sigma}^i, \bar{z}_i) = \textsc{Not Found} \ $}
{
	               $(\bar{\sigma}^i, \bar{z}_i) \gets \textsc{Balas-Simonetti}(\sigma^i_1 \oplus \dots \oplus \sigma^i_{b_i})$\;
	               
		        $\textsc{Store}((\bar{\sigma}^i, \bar{z}_i),H(\sigma^i_1 \oplus \dots \oplus \sigma^i_{b_i}))$\;
}

 $z_\textsc{move} \gets z_\textsc{move} + \bar{z}_i$

}

\vspace*{0.3cm}
\tcp{Filter non-improving moves:}

           \If{$z_\textsc{move} \geq z_\textsc{before}$}
           {\vspace*{0.05cm} \Continue.}

           \vspace*{0.3cm}
           \tcp{At this stage, apply $\phi$ since it is an improving move in $\SB{k}$:}
           Set $x^{t+1} = \phi(x)$ ; $t = t+1$\;
          Replace the routes ($\sigma^1$,$\sigma^2$) by ($\bar{\sigma}^1$,$\bar{\sigma}^2$) in $x^{t+1}$ 
         }
       }
       \Return{$x^t$}
     \end{algorithm}

\noindent
\textbf{Neighborhood Reductions.} First of all, as in the majority of recent local search based metaheuristics for the CVRP, the neighborhood $\mathcal{N}(x^t)$ of each incumbent solution $x^t$ is limited to moves that involve close vertices (Algorithm 1, Line 3). In particular, we use the classical intra-route and inter-route \textsc{Relocate} and \textsc{Swap} moves, for single vertices or generalized to pairs of consecutive vertices, as well as the \textsc{2-Opt} and \textsc{2-Opt*} moves. The resulting neighborhood contains a quadratic number of moves. As detailed in \cite{Vidal2012c}, and in a similar way as \cite{Johnson1997} and \cite{Toth2003}, the search can be restricted to a subset of these moves that reconnect at least one vertex $i$ with a vertex $j$ belonging to the $\Gamma$ closest vertices of $i$. The neighborhood size becomes $\cO(\Gamma n)$, enabling a significant speedup for large-scale problem instances.
\\

\noindent
\textbf{Dynamic move filters.}
To further restrict the search to promising moves, each move $\phi$'s feasibility is evaluated in $\cO(1)$, in terms of capacity constraints, being discarded if it leads to an infeasible solution. The total cost $z(\phi(x^t))$ of the solution generated by $\phi$ prior to its optimization by the $\BS$ decoder is evaluated subsequently. This cost represents an upper bound for the final cost of the move in $\SB{k}$ after the application of the decoder. The move evaluation is pursued only if the solution cost has increased by a factor $1+\psi$ or less due to its application, that is, only if Condition (\ref{condition}) is satisfied. Otherwise, the move is discarded.
\begin{equation}
z(\phi(x^t))  \leq  (1 + \psi) \times z(x^t)
\label{condition}
\end{equation}

Parameter $\psi$ plays an important role in defining how many moves are evaluated.
The higher the value of $\psi$, the less pruning is induced by Equation~(\ref{condition}).
Contrastingly, when $\psi = 0$ only immediately improving neighbors are evaluated.
 Defining a good value for $\psi$ is non-trivial, given that it is an instance-dependent parameter.
 Since a fixed value would not suit instances with different sizes and characteristics, we suggest to use an adaptive parameter.
The principle consists in adjusting $\psi$ to ensure a target range $[\xi^-,\xi^+]$ for the fraction of filtered moves.
 After each 1,000 move evaluations, the fraction $\xi$ of filtered moves is collected and whenever it falls outside of the desired range, $\psi$ is updated.
If this fraction is too large, then $\psi$ is increased by a multiplicative factor $\alpha$. 
Conversely, if $\xi$ is insufficient, then the parameter $\psi$ is decreased:
\begin{equation}
\label{adaptation}
\psi = 
\begin{cases}
\psi \times \alpha & \hspace*{0.3cm} \text{ if } \xi \leq \xi^-, \\
\psi  \ / \ \alpha & \hspace*{0.3cm} \text{ if } \xi \geq \xi^+,\\
\psi  & \hspace*{0.3cm} \text{ otherwise.}
\end{cases}
\end{equation}

\vspace*{0.5cm}
\noindent
\textbf{Global memory.}
The $\BS$ algorithm, used as a decoder, requires a computational effort which grows linearly with the route size and exponentially with parameter $k$. 
It is thus essential to restrain the use of this procedure to a strict minimum and avoid decoding twice the same route over the course of the search. 
To that end, we rely on a global memory to store the routes that have been decoded, avoiding recalculations (Lines 12--14). We use a hashtable for this task, since it allows $\cO(1)$ queries given the hash key associated with a route. 

Two important aspects should be discussed.
First, since the available memory space is finite, some strategy is necessary to limit the memory size in case of an excessive space consumption. To that end, we define an upper bound $\mathcal{M}^\textnormal{max}$ on the number of routes stored in the memory and eliminate half of the entries, those used with less frequency, whenever this limit is attained.

The second aspect to be discussed concerns the effort spent querying the memory.
It is well known that local searches for the CVRP evaluate millions of moves, and that constant-time move evaluations are essential for a good performance.
Capacity checks (Line 5) and simple distance computations (Line 7) can be easily achieved in $\cO(1)$ using incremental move evaluations or concatenation techniques \citep{Vidal2012b}. Moreover, querying the memory for a given route supposes the availability of a hash index which characterizes the associated sequence of visits, but a direct approach that sweeps through the route to compute this index already takes $\cO(n)$ time.
To avoid this bottleneck, we employ specific hash functions and calculation techniques based, again, on concatenations in $\cO(1)$. These concepts are discussed in the following section.

\subsection{Constant-Time Evaluations}
\label{sec:constant}

We use the concatenation strategy of \cite{Vidal2012b,Vidal2015b} to perform efficient cost- and load-feasibility evaluations. This strategy exploits the fact that any route obtained from a classical move $\phi(x^t)$ on an incumbent solution $x^t$ corresponds to a recombination of a bounded number of (customer and depot) visit sequences of $x^t$. As such, the new routes can be expressed as a concatenation of sequences $\sigma_1 \oplus \dots \oplus \sigma_b$.
We also extend this approach to enable $\cO(1)$ computations of hash keys.

To efficiently evaluate the cost, load, and compute the hash keys, we perform a preliminary preprocessing on the $\cO(n^2)$ subsequences of consecutive visits which compose the solution $x^t$.
Four quantities are calculated: the total demand $Q(\sigma)$ of a sequence $\sigma$, its distance $C(\sigma)$, and its hash keys $H^p(\sigma)$ and $H^s(\sigma)$.  For a sequence $\bar{\sigma} = [ i ]$ containing a single visit $i$ with demand~$q_i$, $Q(\bar{\sigma}) = q_i$, $C(\bar{\sigma}) = 0$, $H^p(\bar{\sigma}) = \rho \times i$ and $H^s(\bar{\sigma}) =  \rho^i$, where $\rho$ is a prime number.
Moreover, Equations (\ref{cvrp:Q}--\ref{cvrp:Hp}) extend these quantities, by induction, for any sequence of visits $\sigma_1 \oplus \sigma_2$ expressed as the concatenation of two sequences $\sigma_1$ and~$\sigma_2$. In these equations, $d_{ij}$ expresses the distance between visits $i$ and $j$.
\begin{align}
\label{cvrp:Q}
Q(\sigma^1 \oplus \sigma^2) &= Q(\sigma^1) +  Q(\sigma^2)\\
\label{cvrp:C}
C(\sigma^1 \oplus \sigma^2) &= C(\sigma^1) + d_{\sigma^1(|\sigma^1|), \sigma^2(1)} +  C(\sigma^2) \\
\label{cvrp:Hp}
H^p(\sigma^1 \oplus \sigma^2) &= H^p(\sigma^1) + \rho^{|\sigma_1| } \times H^p(\sigma^2) \\
\label{cvrp:Hs}
H^s(\sigma^1 \oplus \sigma^2) &= H^s(\sigma^1) + H^s(\sigma^2).
\end{align}

As in \cite{Vidal2012b}, Equations (\ref{cvrp:Q}--\ref{cvrp:Hs}) are first employed iteratively, in lexicographic order, to obtain information concerning all sequences during the preprocessing phase. Afterwards, the same equations are used for move evaluations.  Since any route obtained from a classical move corresponds to the concatenation of a bounded number of sequences, it is possible to obtain the associated load, distance, and hash keys by applying these equations a limited number of times. Then, the information on subsequences is updated every time an improving move is applied, a rare occurrence in comparison to the number of moves evaluated.

The two hash functions defined in Equations (\ref{cvrp:Hp}--\ref{cvrp:Hs}) are employed together as a means of reducing chances of two distinct sequences having identical hashes.
The function $H^p$ is a multiplicative hash which depends on the visit permutation \mbox{\citep{Knuth1973}}. Note that, when implementing such a function, the values $\rho^i$ must be precomputed and bounded (taking the rest of the integer division by a large number) to prevent overflow during multiplication. The second function  $H^s$ is an additive hash which only depends on the set of visited customers, and not on the visit sequence.
These hash functions are easily recognized when reformulated as follows:
    \begin{align}
            \label{cvrp:Hp-full}
                H^p(\sigma) &= \sum_{i=1}^{|\sigma|} \rho^i \times \sigma_{i} \\
            \label{cvrp:Hs-full}
                H^s(\sigma) &= \sum_{i=1}^{|\sigma|} \rho^{\sigma_{i}}.
        \end{align}

These functions fit well our purposes due to their inductive definition based on the concatenation operation. They are employed, along with the route distance and its number of visits, to verify a correct match in the memory in $\cO(1)$ without a complete route comparison in $\cO(n)$. To minimize the risk of two routes having identical hashes, we duplicated these hash functions with different values for $\rho$, leading to four hash values overall.
In the first case, $\rho$ is set to the smallest prime number greater than the number of customers. In the second case, $\rho$ is set to $31$ \citep[multiplier of][]{Kernighan1988}.
Despite this strategy, a tiny chance of false positives remains. However, no false positive was registered within our computational experiments considering multiple runs on 100 different instances.

\subsection{Reshaping the search space -- Tunneling strategy}
\label{cvrp:sec:reshaping}

Until now, the purpose of the global memory has been focused on saving computational effort by avoiding duplicate calls to the $\BS$ algorithm. Yet, as shown in the following, this structure can be exploited to a larger extent to promote the discovery of good solutions.

Consider two input routes $\sigma_1$ and $\sigma_2$ indexed in the order of their appearance in the local search, and representing the same set of customer visits.
The route $\sigma_1$ is first saved in the global memory along with its associated \Bk--optimal route of cost $z_1$. Subsequently, the LS considers $\sigma_2$ without finding a match in the memory, triggering a new execution of the $\BS$ algorithm and leading to a cost $z_2$.
if $z_2 > z_1$, then the algorithm has failed to recognize that a better TSP tour has been found in prior search for the same customers.

\begin{figure}[htb]
\centering
\vspace*{0.2cm}
\includegraphics[width=0.65\textwidth]{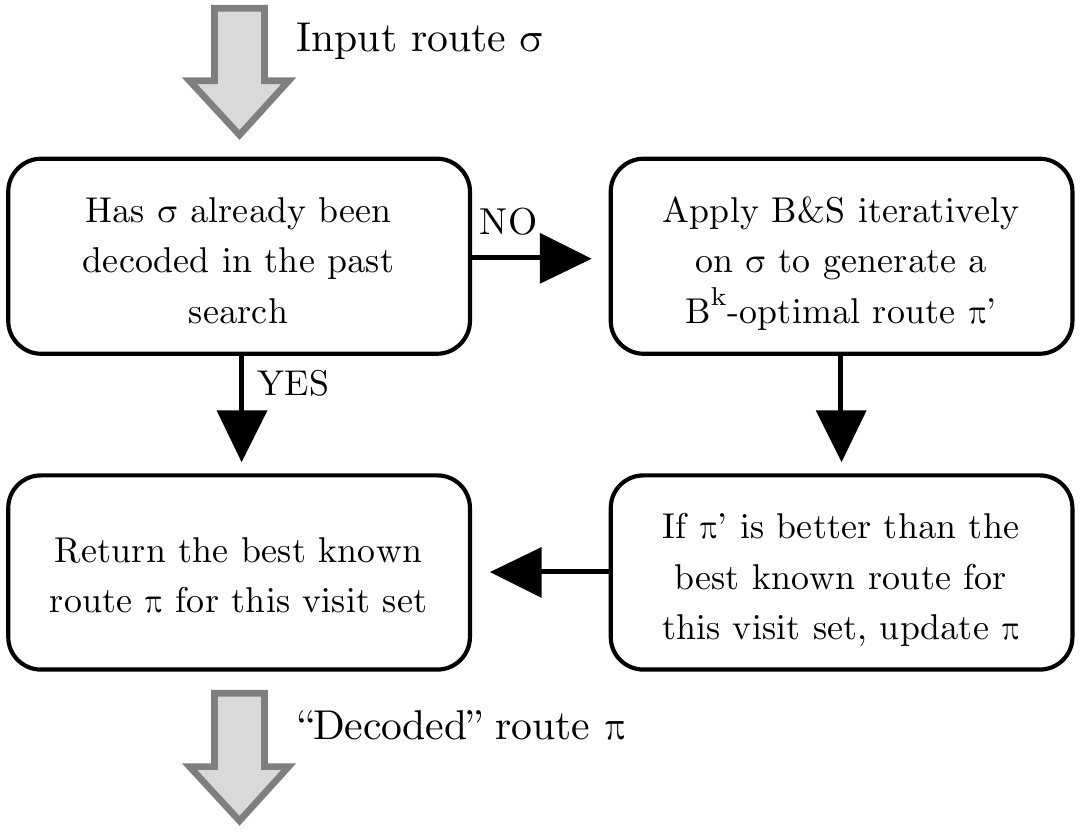}
\caption{Tunneling strategy and solution decoding}
\label{fig:caching_strategies}
\end{figure}

To improve the behavior of the algorithm in such situations, we introduce a guidance mechanism called \emph{tunneling}. Guidance techniques are a set of strategies which analyze and exploit the search history to direct the search towards promising or unexplored regions of the search space \citep{Crainic2008}. Our method works as follows: every route $\sigma$ issued from a non-filtered LS move and absent from the memory is decoded by the $\BS$ algorithm; yet, rather than directly returning the output of $\BS$, the algorithm finds and returns the best known permutation of the visits for this customer set, found over previous $\BS$ executions. The goal of this strategy is to intensify the search around known high-quality tours without jeopardizing the discovery of better route configurations. It can be efficiently implemented with a refinement of the hashtable-based memory structure, by grouping the routes into different buckets according to their visit set, and using the additive hash function of Equation~(\ref{cvrp:Hs}) for~$\cO(1)$ queries. Figure \ref{fig:caching_strategies} summarizes this process.

This tunneling strategy has a significant impact on the search space.
Initially, as the search starts, the algorithm explores the space $\SB{k}$.
Then, as the search progresses, the memory starts to be filled, and the algorithm re-introduces more and more frequently the best known routes in its solutions. In a hypothetical situation where all feasible routes have already been memorized (hypothetical due to the needed exponential memory size), the TSP--optimal routes would be systematically returned, and the algorithm behaves as if it was searching in $\SA$. With this limit case in mind, the tunneling strategy contributes to reshape the space $\SB{k}$ into $\SA$ as the search progresses. Moreover, note that this strategy remains fully relevant for a classical neighborhood search, even without $\BS$ decoder.

\begin{figure}[!ht]
\centering
\includegraphics[width=0.85\textwidth,page=4]{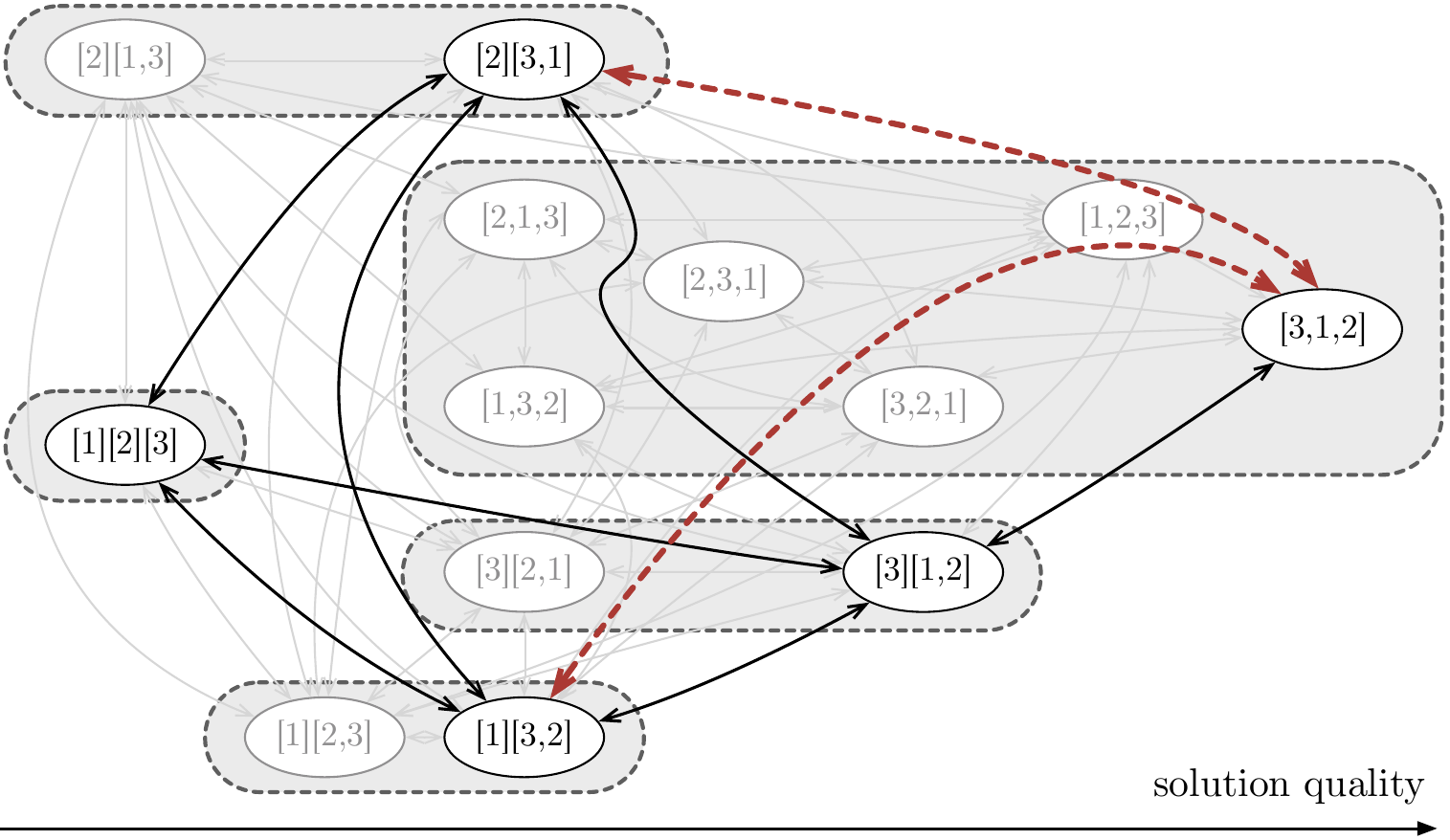}
\caption{Search space $\cS_1^{\textsc{B}}$ with the tunneling search strategy, after discovering the route $[3,1,2]$}
\label{fig:tunneling}
\end{figure}

Consider the same example as depicted in Section \ref{sec:search_space} (Figures \ref{fig:neigh}, \ref{fig:neigh_concorde} and \ref{fig:neigh_bs}).
Suppose that the route $[3,2,1]$ has been identified in the past search along with its associated \Bk--optimal tour $[3,1,2]$ and that the tunneling strategy is employed. Figure \ref{fig:tunneling} illustrates the resulting search space: all solutions that include `$[1,2,3]$' in their neighborhood now point towards solution  `$[3,1,2]$' instead, which has identical customer-to-vehicle \textsc{Assignments} but a lower cost due to better \textsc{Sequencing} decisions. With only one route in memory, the resulting search space already becomes equivalent to $\SA$.

\section{Computational experiments}
\label{sec:experiments}

We conducted extensive computational analyses to measure the benefits of a search in $\SA$~and~$\SB{k}$, the impact of the tunneling strategy and move evaluation filters. For these tests, we considered the simple local search (LS) described in Section~\ref{sec:methodology}, as well as a more advanced metaheuristic, the UHGS of \cite{Vidal2012,Vidal2012b}, which was adapted by replacing its native local search by the proposed method. This extension of the method will be referred as UHGS-BS. Except this modification of the local search, all other parameters and procedures of UHGS-BS remain the same as in the original article, and the termination criterion is set to $It_\textsc{max} = 20000$ consecutive iterations without improvement.

A single local search requires only a limited computational effort. Thus, the first analyses (Section \ref{tests1}) on the performance of the LS in $\SA$ and $\SB{k}$ could be done for multiple values of the range parameter $k$. It was also possible to evaluate the impact of $k$ without the interference of dynamic move filters and tunneling techniques.
Since UHGS-BS performs a more extensive search with multiple local search runs, our analyses with this method (Section \ref{tests2}) are focused on a smaller set of values for the range parameter $k$, in the presence of dynamic move filters. 

All algorithms were implemented in C++ and executed on a single thread of an Intel(R) Xeon(R) E5-2680v3 CPU. A RAM limit of 8GB was imposed for each run. To solve the TSP problems when considering the space $\SA$, we used the \textsc{Concorde} solver \citep{Concorde2006}. To conduct the experiments with the UHGS-BS, we used the code base made available at \url{https://github.com/vidalthi/HGS-CARP}, from \cite{Vidal2015a}.

\subsection{Preliminary experiments with a simple local search}
\label{tests1}

In a first experiment, we tested the local search of Section~\ref{sec:methodology} on spaces $\SA$ and $\SB{k}$ for $k \in \{0,\dots,9\}$, to observe the growth of its computational time as a function of the parameter~$k$, and identify a range of values over which the approach remains practical. As an initial solution, we used the result of the \emph{savings} algorithm of \cite{Clarke1964}. We set $\psi=\infty$ in order to observe the results without the interference of move filters.

\begin{figure}[!h]
\centering
\includegraphics[width=1.0\textwidth,page=1]{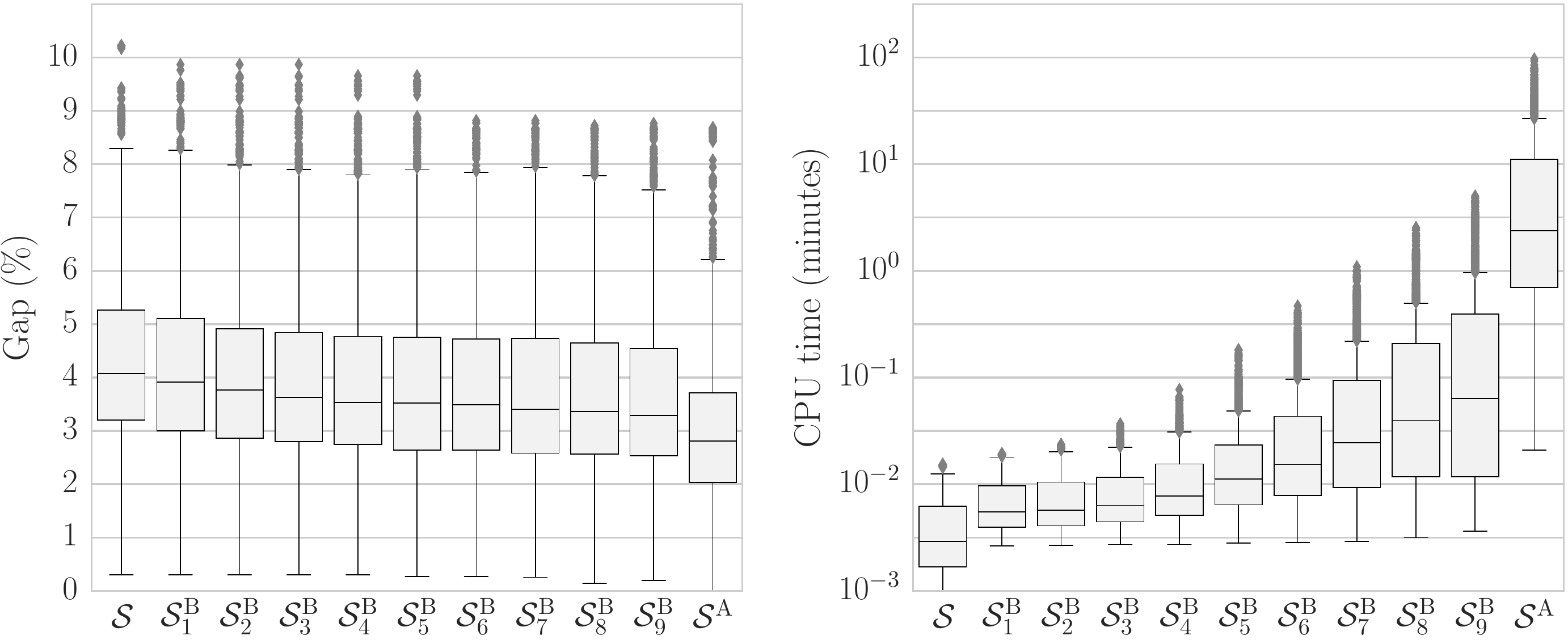}
\caption{Solution quality and CPU time of the LS, depending on the search space}
\label{cvrp:boxplot:ls}
\end{figure}

We considered the 100 recent benchmark instances of \cite{Uchoa2017}, as these instances remain highly challenging for metaheuristics and cover a larger variety of instance size and characteristics: demand and customer distribution, depot location, and route length. For each instance, we ran the method 20 times with different random seeds. The results are summarized in Figure~\ref{cvrp:boxplot:ls}, in the form of boxplots. The leftmost graph represents the percentage gap in terms of solution quality, relative to that of the best known solution (BKS) collected from \cite{Uchoa2017}: $\text{Gap} = 100 \times (z - z_\textsc{bks}) / z_\textsc{bks}$, where $z$ is the solution value of the method and $z_\textsc{bks}$ is the BKS value. The rightmost graph represents the CPU time of the method, using a logarithmic scale.

As illustrated by these experiments, the search in space $\SA$ visibly leads to solutions of higher average quality, albeit in a CPU time largely greater than that of a classical local search in $\cS$. The solution quality of a search in the space $\SB{k}$ consistently increases as $k$ grows, along with the needed CPU time. When $k=0$, the algorithm behaves as a classical local search in $\cS$. When $k$ is large, the method becomes more similar to a search in $\SA$. A difference of solution quality can still be noticed between $\SB{9}$ and $\SA$, due to some instances containing up to 25 deliveries per route.

Each value of the range parameter $k$ establishes a trade-off between computational effort and solution quality. The computational time of the local search does not exceed two seconds when exploring $\SB{k}$ with $k \leq 2$, but it culminates to 1000 seconds on the largest instances when exploring $\SA$. Clearly, the additional CPU time required by Concorde is not compatible with the repeated application of local searches within state-of-the-art metaheuristic algorithms, such that we should concentrate our efforts on the exploration of the space $\SB{k}$ with moderate values of~$k$.

\begin{figure}[!ht]
\vspace{10pt}
(a) Instances with average route cardinality in range [3.0, 4.55]:\vspace{5pt}\\
\begin{subfigure}{1.0\textwidth}
  \centering
  \includegraphics[width=1.0\textwidth,page=1]{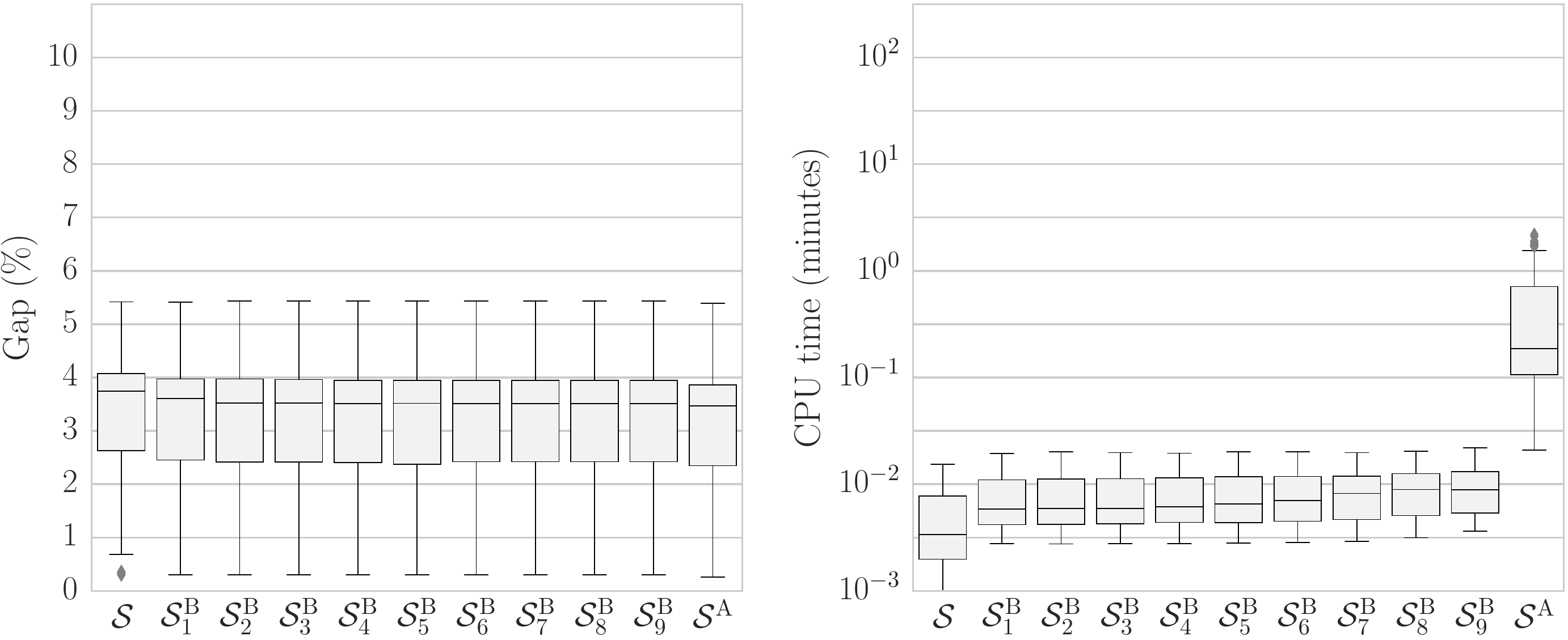}
\end{subfigure}
\vspace{10pt}\\
(b) Instances with average route cardinality in range [16.47, 24.43]:\vspace{5pt}\\
\begin{subfigure}{1.0\textwidth}
  \centering
  \includegraphics[width=1.0\textwidth,page=1]{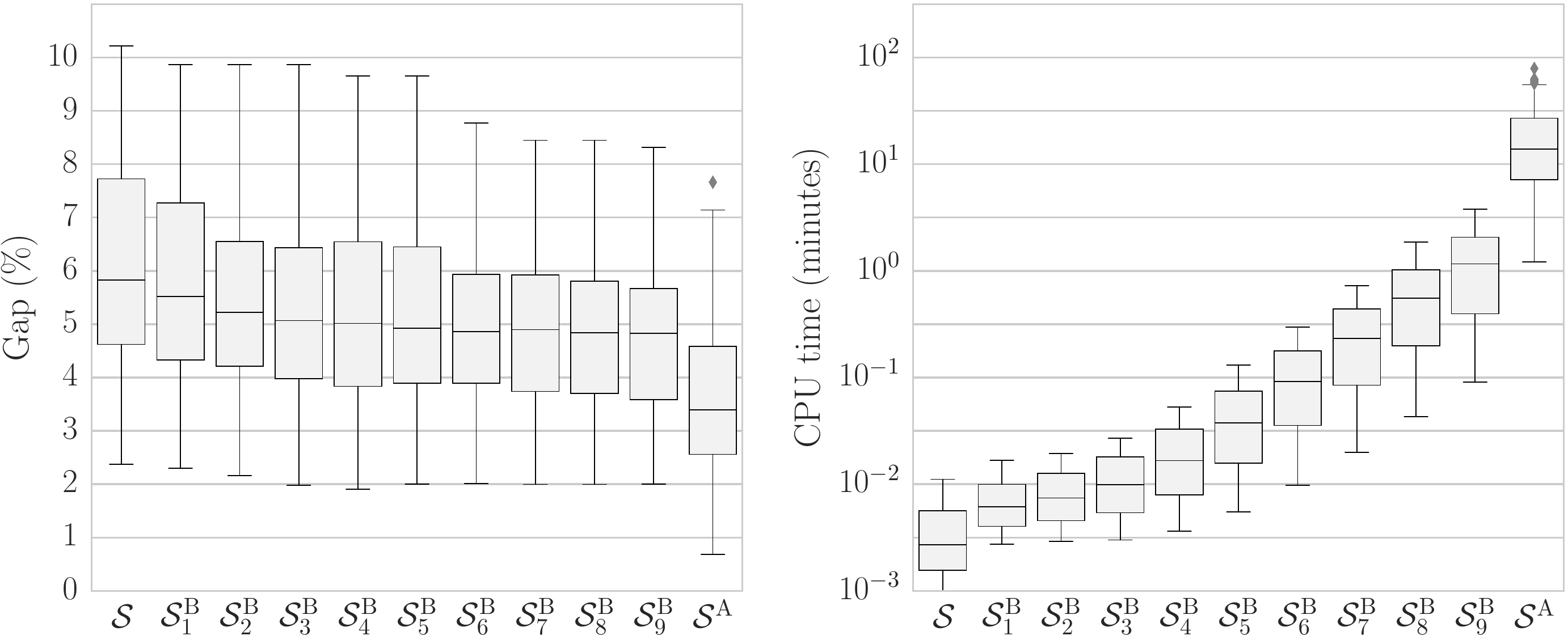}
\end{subfigure}
\caption{Results (solution quality and runtime) of local search on different search spaces for instances with different route cardinalities}
\label{cvrp:boxplot:ls_groups}
\end{figure}

Moreover, a careful analysis of these results on subsets of instances with a different average number of customers per route gives additional insights. This analysis is reported in Figure~\ref{cvrp:boxplot:ls_groups}, considering the 20 instances with smallest average route cardinality, in the range $[3.0,4.55]$, and the 20 instances with largest average route cardinality, in $[16.47,24.43]$.

As illustrated in Figure~\ref{cvrp:boxplot:ls_groups}, the benefits of a search in $\SB{k}$ or $\SA$ are small for instances with a small number of customer visits per route. In particular, all runs with $k \geq 4$ lead to a similar solution quality and CPU time. This is due to the fact that $\BS$ does an exact TSP optimization when $k$ is greater or equal to the route cardinality. In contrast, the benefits in terms of solution quality are larger on instances with a high number of customer visits per vehicle. We observe a significant improvement of the solutions when $k$ varies from $0$ to $4$, from an average gap of $6.15\%$ down to $5.25\%$.
Subsequently, as $k$ increases beyond $4$, the rate of improvement is smaller. Increasing $k$ up to the maximum route size would still be beneficial, but impracticable in terms of CPU time.

\subsection{Experiments with UHGS-BS -- range parameter, move filters and tunneling}
\label{tests2}

As viewed in the previous section, a local search in the space $\SB{k}$ can lead to solutions of better quality than a search in $\cS$, at the expense of a higher computational effort. Still, even if solution improvements were observed for simple local searches, it is an open question whether the inclusion of these extended search procedures into state-of-the-art metaheuristics translates back into significant quality improvements.

This section now analyses the performance of the UHGS-BS metaheuristic, considering $k \in \{1,\dots,5\}$ in combination with dynamic move filters. The move filters are designed to eliminate a large fraction of complete move evaluations, such that a larger computational effort can be spent in the evaluation of each remaining move without significantly impacting the overall CPU time of the method.

After preliminary experiments, we observed that setting $k=2$ and $[\xi^-,\xi^+] = [90\%,95\%]$ establishes a good compromise between solution quality and computational effort. This configuration, without tunneling strategy, was used as baseline for the experiments of this section. We then modified each parameter and design choice, in turn, to investigate its impact. This experimentation was restricted to a subset of 40 instances with a number of customer visits $n \in \{195,411\}$, as these instances require limited CPU time and remain challenging for state-of-the-art metaheuristics. For each instance, ten runs have been conducted with different~seeds.\\

\begin{figure}[!h]
\centering
\includegraphics[width=1.0\textwidth,page=1]{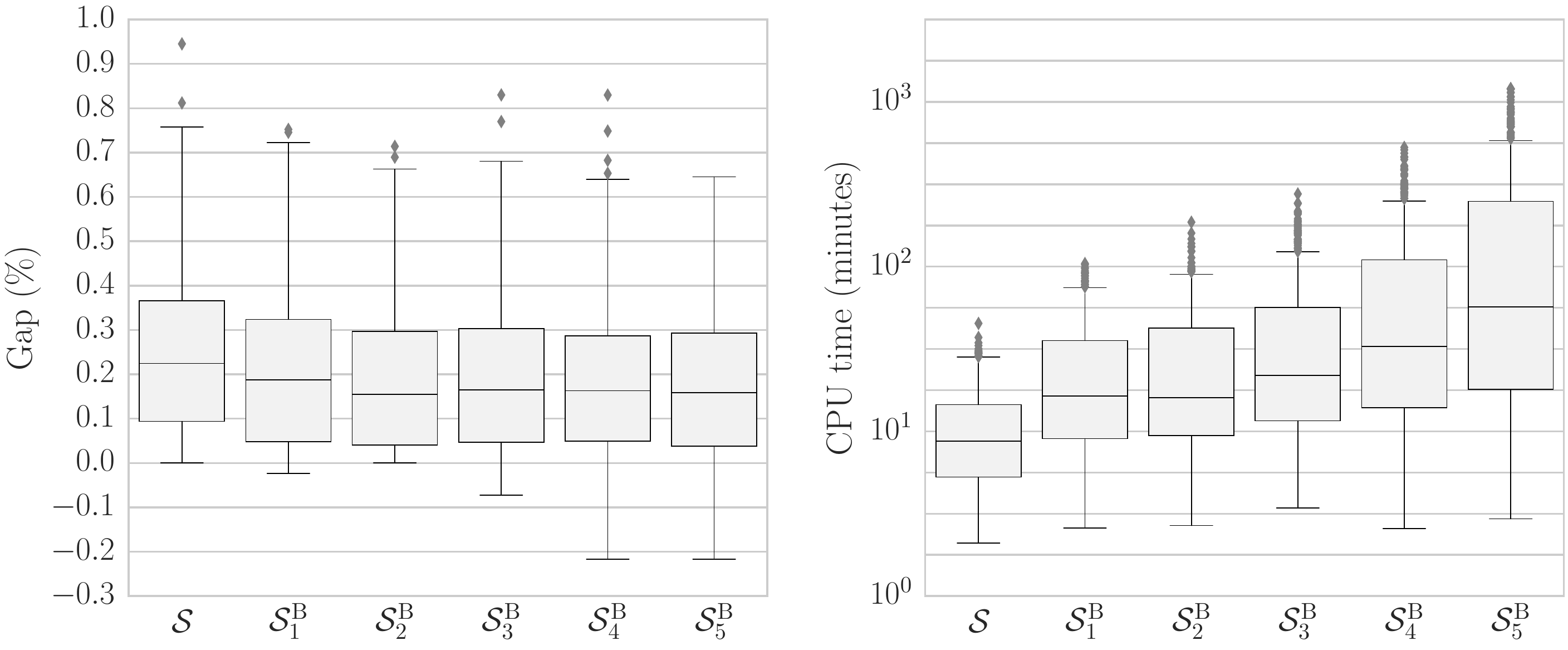}
\caption{Results (solution quality and runtime) of UHGS-BS for different cache strategies and $k$ values}
\label{cvrp:boxplot:k}
\end{figure}

\noindent
\textbf{Impact of the range parameter.}  Figure~\ref{cvrp:boxplot:k} compares the performance, in terms of average percentage gap and CPU time, of the classical UHGS with that of its extended version searching in $\SB{k}$, for $k \in \{1,\dots,5\}$.
As in previous figures, the results are presented in the form of boxplots. 

As expected, the gaps obtained with the variants of UHGS-BS are much smaller than those of simple local searches, due to the better exploration capabilities of the method. Average gaps range from $0.25\%$ when exploring the classical search space $\cS$, to $0.18\%$ when exploring $\SB{5}$. As highlighted by pairwise Wilcoxon tests, statistically significant differences exist between the results of $\cS$, $\SB{1}$ and $\SB{2}$ (p-values $< 0.05$). A \emph{decreasing returns} effect can also be observed; the difference of quality between the solutions obtained by $\cS$ and $\SB{1}$ is larger than between $\SB{1}$ by $\SB{2}$, which is turn is larger than between $\SB{2}$ by $\SB{3}$, and so on.
The configuration $\SB{3}$, in particular, achieves a good trade-off between quality and search effort. As demonstrated by the outliers with negative gap in the figure, this configuration has led to new best known solutions for surprisingly small instances with 256 and 294 customers, on which hundreds of test runs had been conducted in past. This is an indication that the search in $\SB{k}$ has the potential to lead to structurally different solutions, which were not attained with more classical searches. The main challenge therefore is to harness this ability without sacrificing too much computational effort.\\

\begin{figure}[!h]
\centering
\includegraphics[width=1.0\textwidth,page=1]{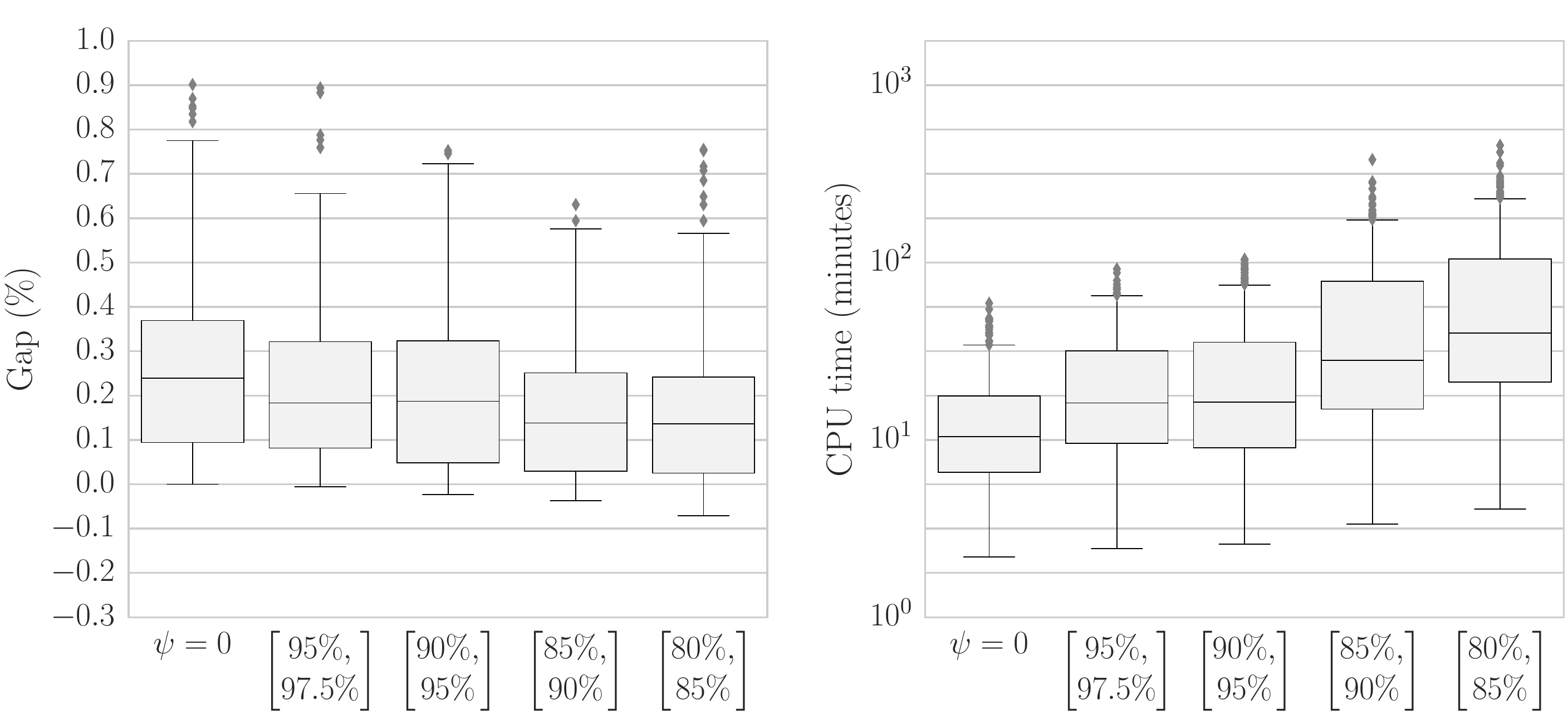}
\caption{Results (solution quality and runtime) of UHGS-BS for different $[\xi^-,\xi^+]$ values}
\label{cvrp:boxplot:lb}                 
\end{figure}

\noindent
 \textbf{Impact of the dynamic move filter.} 
Figure~\ref{cvrp:boxplot:lb} investigates the impact of different target intervals $[\xi^-,\xi^+]$ (desired quantity of filtered moves -- Section \ref{sec:lower_bound}) for the dynamic move filters. The range parameter remains fixed to $k = 2$. It also indicates the results obtained when filtering \emph{all} non-improving moves ($\psi=0$), which is equivalent to using $\BS$ only as a post-optimization procedure, after the discovery of each improving move.

These experiments demonstrate that move filters have a large incidence on the solution quality. These filters, are however, essential to maintain a low computational effort. In particular, filtering all non-improving moves prior to the evaluation of {\BS} ($\psi=0$) leads to an average gap of 0.26\%, compared to 0.21\% when setting $[\xi^-,\xi^+] = [90\%,95\%]$ as a target for the dynamic move filter and evaluating the non-filtered moves in combination with {\BS}.
This validates an important hypothesis explored in this article: many moves that are usually discarded in regular local search methods and metaheuristics can lead to improved solutions when applied in combination with a route optimization procedure. Naturally, this capability goes along with an increased CPU time. Nonetheless, by an adequate calibration of the move filters, the total CPU time dedicated to {\BS} can be restricted enough to not intervene as a bottleneck. This is visible by the results of configuration  $[\xi^-,\xi^+] = [90\%,95\%]$, which on average used no more than twice the time of $\psi=0$. For the remainder of these analyses, we selected this configuration, which establishes a good balance between the exploitation of the capabilities of {\BS} and the computational effort.\\

\noindent
 \textbf{Impact of the tunneling strategy.}
Finally, Table~\ref{table_tunneling} compares the performance of UHGS-BS without and with the tunneling strategy. The range parameter has been set to $k=2$, and $[\xi^-,\xi^+] = [90\%,95\%]$.
The leftmost group of columns reports the instance name, number of customers $n$ average number of visits per route in the BKS for each instance. Then, the next columns present, for each approach, the average solution value over ten runs, best solution value, average CPU time, percentage of routes which have been successfully queried from the global memory without a re-evaluation, number of executions of the {\BS} optimization procedure, total number of iterations.

\begin{table}[htbp]\scriptsize
  \centering
  \caption{Impact of the tunneling strategy for 40 medium-size instances}
  \label{table_tunneling}
  \def\arraystretch{1.05}\tabcolsep=2.6pt
\begin{tabular}{rrrrrrrrrrrrrrrr}
\toprule
\multirow{2}[4]{*}{\#} & \multicolumn{1}{l}{\multirow{2}[4]{*}{Instance}} &       & \multicolumn{6}{c}{Without tunneling strategy} &       & \multicolumn{6}{c}{With tunneling strategy} \\
\cmidrule{4-9}\cmidrule{11-16}      &       &       & Time  & Avg.  & Best  & Cache & \BS & Iters &       & Time  & Avg.  & Best  & Cache & \BS & Iters \\
\midrule
21    & \multicolumn{1}{l}{X-n195-k51} &       & 3.8   & 44289.5 & \textbf{44225} & 60.6\% & 1.4E+08 & 34577 &       & 5.2   & \textbf{44272.6} & \textbf{44225} & 60.2\% & 1.8E+08 & 41713 \\
22    & \multicolumn{1}{l}{X-n200-k36} &       & 5.3   & \textbf{58618.5} & \textbf{58578} & 97.8\% & 1.1E+07 & 43664 &       & 6.2   & 58625.7 & \textbf{58578} & 97.8\% & 1.0E+07 & 41788 \\
23    & \multicolumn{1}{l}{X-n204-k19} &       & 5.8   & 19569.0 & \textbf{19565} & 91.0\% & 2.7E+07 & 31315 &       & 6.6   & \textbf{19568.5} & \textbf{19565} & 92.6\% & 2.2E+07 & 31118 \\
24    & \multicolumn{1}{l}{X-n209-k16} &       & 15.5  & \textbf{30673.1} & \textbf{30656} & 77.6\% & 1.1E+08 & 35319 &       & 17.8  & 30676.8 & \textbf{30656} & 81.0\% & 9.1E+07 & 36174 \\
25    & \multicolumn{1}{l}{X-n214-k11} &       & 34.8  & 10873.2 & \textbf{10856} & 70.5\% & 2.1E+08 & 48054 &       & 32.3  & \textbf{10872.0} & \textbf{10856} & 76.2\% & 1.5E+08 & 46133 \\
26    & \multicolumn{1}{l}{X-n219-k73} &       & 12.9  & \textbf{117602.7} & \textbf{117595} & 3.7\% & 5.0E+08 & 30949 &       & 11.0  & 117603.8 & \textbf{117595} & 3.6\% & 5.0E+08 & 30773 \\
27    & \multicolumn{1}{l}{X-n223-k34} &       & 6.2   & 40490.3 & \textbf{40437} & 95.4\% & 2.3E+07 & 38698 &       & 7.7   & \textbf{40489.0} & \textbf{40437} & 95.7\% & 2.1E+07 & 39268 \\
28    & \multicolumn{1}{l}{X-n228-k23} &       & 14.5  & 25784.1 & \textbf{25743} & 87.9\% & 8.5E+07 & 48303 &       & 14.6  & \textbf{25779.5} & \textbf{25743} & 88.6\% & 6.9E+07 & 42435 \\
29    & \multicolumn{1}{l}{X-n233-k16} &       & 10.9  & \textbf{19305.3} & \textbf{19230} & 87.1\% & 4.7E+07 & 31274 &       & 9.5   & 19309.0 & \textbf{19230} & 90.4\% & 3.2E+07 & 29652 \\
30    & \multicolumn{1}{l}{X-n237-k14} &       & 16.6  & 27053.8 & 27050 & 75.9\% & 9.7E+07 & 26927 &       & 17.5  & \textbf{27047.2} & \textbf{27042} & 81.1\% & 8.1E+07 & 30294 \\
31    & \multicolumn{1}{l}{X-n242-k48} &       & 8.7   & 82922.0 & \textbf{82792} & 94.3\% & 4.4E+07 & 49734 &       & 12.9  & \textbf{82919.6} & \textbf{82792} & 93.6\% & 6.6E+07 & 62591 \\
32    & \multicolumn{1}{l}{X-n247-k47} &       & 14.6  & \textbf{37380.7} & \textbf{37278} & 49.2\% & 6.1E+08 & 63084 &       & 13.7  & 37393.1 & \textbf{37281} & 49.1\% & 5.1E+08 & 52872 \\
33    & \multicolumn{1}{l}{X-n251-k28} &       & 12.7  & \textbf{38767.4} & \textbf{38699} & 90.8\% & 6.4E+07 & 48107 &       & 11.4  & 38793.5 & \textbf{38699} & 91.2\% & 4.5E+07 & 35265 \\
34    & \multicolumn{1}{l}{X-n256-k16} &       & 8.8   & 18880.0 & 18880 & 85.9\% & 4.3E+07 & 22991 &       & 9.4   & \textbf{18875.9} & \textbf{18839} & 89.6\% & 3.3E+07 & 24930 \\
35    & \multicolumn{1}{l}{X-n261-k13} &       & 43.0  & 26628.7 & \textbf{26579} & 71.9\% & 2.3E+08 & 46952 &       & 38.9  & \textbf{26620.5} & \textbf{26586} & 75.1\% & 1.7E+08 & 39770 \\
36    & \multicolumn{1}{l}{X-n266-k58} &       & 14.3  & \textbf{75703.6} & \textbf{75558} & 98.9\% & 1.5E+07 & 69013 &       & 17.8  & 75737.9 & \textbf{75646} & 98.9\% & 1.5E+07 & 72435 \\
37    & \multicolumn{1}{l}{X-n270-k35} &       & 6.6   & \textbf{35310.5} & \textbf{35303} & 96.8\% & 1.5E+07 & 34796 &       & 8.5   & 35318.5 & \textbf{35303} & 97.1\% & 1.4E+07 & 35654 \\
38    & \multicolumn{1}{l}{X-n275-k28} &       & 11.7  & 21256.9 & \textbf{21245} & 89.0\% & 5.9E+07 & 36470 &       & 12.8  & \textbf{21255.0} & \textbf{21245} & 91.0\% & 4.4E+07 & 34488 \\
39    & \multicolumn{1}{l}{X-n280-k17} &       & 69.1  & 33614.0 & 33593 & 66.7\% & 4.1E+08 & 53528 &       & 86.3  & \textbf{33600.8} & \textbf{33584} & 70.1\% & 4.3E+08 & 63405 \\
40    & \multicolumn{1}{l}{X-n284-k15} &       & 82.2  & 20306.0 & 20255 & 70.0\% & 4.3E+08 & 73605 &       & 67.4  & \textbf{20283.3} & \textbf{20238} & 73.9\% & 3.0E+08 & 61772 \\
41    & \multicolumn{1}{l}{X-n289-k60} &       & 18.3  & 95501.4 & 95395 & 69.6\% & 4.5E+08 & 75465 &       & 20.7  & \textbf{95475.6} & \textbf{95245} & 69.5\% & 4.4E+08 & 72940 \\
42    & \multicolumn{1}{l}{X-n294-k50} &       & 9.3   & \textbf{47272.9} & 47240 & 96.8\% & 2.6E+07 & 47162 &       & 11.6  & 47289.7 & \textbf{47239} & 96.3\% & 3.3E+07 & 53125 \\
43    & \multicolumn{1}{l}{X-n298-k31} &       & 7.1   & 34282.8 & \textbf{34231} & 94.6\% & 2.3E+07 & 27625 &       & 8.1   & \textbf{34280.8} & \textbf{34231} & 95.2\% & 2.0E+07 & 26203 \\
44    & \multicolumn{1}{l}{X-n303-k21} &       & 22.3  & 21839.1 & \textbf{21744} & 86.7\% & 1.0E+08 & 43974 &       & 21.1  & \textbf{21833.3} & \textbf{21744} & 88.9\% & 7.9E+07 & 39415 \\
45    & \multicolumn{1}{l}{X-n308-k13} &       & 73.6  & \textbf{25893.2} & \textbf{25864} & 68.5\% & 3.3E+08 & 48194 &       & 63.6  & 25911.8 & \textbf{25866} & 71.8\% & 2.6E+08 & 42804 \\
46    & \multicolumn{1}{l}{X-n313-k71} &       & 15.8  & \textbf{94270.5} & \textbf{94169} & 57.8\% & 4.7E+08 & 54350 &       & 16.8  & 94289.7 & \textbf{94192} & 57.6\% & 4.9E+08 & 58766 \\
47    & \multicolumn{1}{l}{X-n317-k53} &       & 29.7  & \textbf{78389.3} & \textbf{78372} & 97.5\% & 3.3E+07 & 60402 &       & 31.6  & 78405.0 & \textbf{78380} & 97.7\% & 3.2E+07 & 65092 \\
48    & \multicolumn{1}{l}{X-n322-k28} &       & 10.9  & 29918.8 & 29880 & 93.3\% & 3.7E+07 & 36519 &       & 17.4  & \textbf{29894.6} & \textbf{29834} & 93.9\% & 4.3E+07 & 48316 \\
49    & \multicolumn{1}{l}{X-n327-k20} &       & 31.2  & 27618.2 & \textbf{27560} & 81.3\% & 1.5E+08 & 44348 &       & 35.5  & \textbf{27591.7} & \textbf{27560} & 83.4\% & 1.5E+08 & 47634 \\
50    & \multicolumn{1}{l}{X-n331-k15} &       & 71.6  & 31138.6 & \textbf{31103} & 65.2\% & 3.7E+08 & 47298 &       & 63.2  & \textbf{31128.3} & \textbf{31103} & 67.3\% & 3.0E+08 & 41354 \\
51    & \multicolumn{1}{l}{X-n336-k84} &       & 34.3  & \textbf{139522.2} & \textbf{139303} & 39.0\% & 1.4E+09 & 86713 &       & 40.9  & 139572.7 & \textbf{139339} & 39.7\% & 1.4E+09 & 89964 \\
52    & \multicolumn{1}{l}{X-n344-k43} &       & 11.7  & 42161.9 & 42086 & 97.1\% & 2.3E+07 & 42768 &       & 22.7  & \textbf{42136.7} & \textbf{42066} & 97.4\% & 3.0E+07 & 65989 \\
53    & \multicolumn{1}{l}{X-n351-k40} &       & 26.1  & 26010.1 & 25958 & 94.6\% & 7.9E+07 & 72726 &       & 32.8  & \textbf{25991.2} & \textbf{25947} & 94.7\% & 8.3E+07 & 76997 \\
54    & \multicolumn{1}{l}{X-n359-k29} &       & 65.3  & \textbf{51676.6} & 51608 & 84.9\% & 3.1E+08 & 84963 &       & 56.1  & 51684.5 & \textbf{51569} & 85.9\% & 2.4E+08 & 67011 \\
55    & \multicolumn{1}{l}{X-n367-k17} &       & 92.0  & \textbf{22870.1} & \textbf{22814} & 72.2\% & 4.0E+08 & 60550 &       & 66.8  & 22894.6 & \textbf{22814} & 71.1\% & 2.8E+08 & 39572 \\
56    & \multicolumn{1}{l}{X-n376-k94} &       & 61.3  & \textbf{147736.8} & \textbf{147718} & 99.2\% & 1.6E+07 & 64604 &       & 53.1  & 147738.1 & \textbf{147718} & 99.2\% & 1.4E+07 & 61378 \\
57    & \multicolumn{1}{l}{X-n384-k52} &       & 26.0  & 66254.7 & \textbf{66095} & 97.7\% & 4.2E+07 & 77380 &       & 29.1  & \textbf{66184.1} & \textbf{66133} & 97.6\% & 3.8E+07 & 68679 \\
58    & \multicolumn{1}{l}{X-n393-k38} &       & 22.7  & 38361.0 & 38269 & 92.7\% & 7.2E+07 & 47398 &       & 26.9  & \textbf{38305.5} & \textbf{38260} & 93.7\% & 6.5E+07 & 48965 \\
59    & \multicolumn{1}{l}{X-n401-k29} &       & 79.6  & 66357.8 & 66269 & 83.9\% & 3.4E+08 & 86787 &       & 69.7  & \textbf{66347.6} & \textbf{66225} & 86.1\% & 2.7E+08 & 78713 \\
60    & \multicolumn{1}{l}{X-n411-k19} &       & 110.9 & 19753.4 & 19725 & 72.8\% & 4.9E+08 & 65729 &       & 74.8  & \textbf{19748.3} & \textbf{19719} & 77.9\% & 2.9E+08 & 51034 \\
\midrule
\multicolumn{3}{r}{Average values:} & 30.4  & 0.17\% & 0.02\% & 80.2\% & 2.1E+08 & 51058 &       & \textbf{29.3} & \textbf{0.15\%} & \textbf{0.01\%} & 81.5\% & 1.8E+08 & 49912 \\
\bottomrule
\end{tabular}%

\end{table}

The tunneling strategy leads to an average Gap(\%) of $0.15\%$, compared to $0.17\%$ without tunneling. Although this is only a difference of $0.02\%$, reducing the gap indeed becomes harder as the solution quality approaches known BKS and optimal solutions. This improvement in solution quality also comes with a reduction of the overall CPU time, as the tunneling stimulates a faster convergence towards local minima and allows a better management of the global memory, by storing at most one route per customer set. As a consequence, the chances of successful queries in the memory is sensibly higher (81.5\% compared to 80.2\% on average), thereby reducing on average the number of \BS\ executions as well as the total number of iterations. Based on these observations, tunneling is beneficial without any other visible counterpart. We will use this mechanism for the final tests on the complete set of benchmark instances, in the next section.

\FloatBarrier
\subsection{Comparison with recent state-of-the-art algorithms}
\label{tests3}

Finally, this section reports detailed results of UHGS-BS, using the baseline configuration and the tunneling strategy, on the complete set of 100 instances proposed by \cite{Uchoa2017}.
The results of UHGS-BS are compared to that of the current state-of-the-art algorithms: the Hybrid Iterated Local Search (HILS) proposed by \cite{Subramanian2013}, and the original UHGS of \cite{Vidal2012,Vidal2012b}, which were executed 50 times for each instance. To keep the total computational effort within reasonable limits, UHGS-BS was executed 10 times for each instance. The maximum number of consecutive iterations without improvements was set to $It_{\textsc{max}} = 50000$ to evaluate UHGS-BS in the same conditions as UHGS on this set of instances \citep{Uchoa2017}. A hard runtime limit of 24 hours was imposed for each run.

Tables~\ref{cvrp:tab:small} and \ref{cvrp:tab:large} report the results obtained with UHGS-BS, in comparison with UHGS and HILS. The columns present the average CPU time in minutes, the average solution value and the best solution value for all approaches. The best results are highlighted in boldface in the table, with $\circledast$ indicating an improvement over the best known solution from \url{http://vrp.atd-lab.inf.puc-rio.br/}, and the average gap to the best solution produced by the considered methods (including those generated during this research) is presented for each algorithm in each table's final row.

    \begin{table}[htbp]\scriptsize
      \centering
      \caption{Results for the instances with up to 331 customers from \cite{Uchoa2017}}
      \label{cvrp:tab:small}
      \def\arraystretch{1.05}\tabcolsep=3.0pt
\begin{tabular}{rrrrrrrrrrrrrrrr}
\toprule
\multirow{2}[4]{*}{\#} & \multicolumn{1}{l}{\multirow{2}[4]{*}{Instance}} &       & \multicolumn{3}{c}{ILS} &       & \multicolumn{3}{c}{UHGS} &       & \multicolumn{4}{c}{UHGS-BS} \\
\cmidrule{4-6}\cmidrule{8-10}\cmidrule{12-15}      &       &       & Time  & Average & Best  &       & Time  & Average & Best  &       & Time  & Average & Best &  \\
\midrule
1     & \multicolumn{1}{l}{X-n101-k25} &       & 0.1   & \textbf{27591.0} & \textbf{27591} &       & 1.4   & \textbf{27591.0} & \textbf{27591} &       & 2.4   & \textbf{27591.0} & \textbf{27591} &  \\
2     & \multicolumn{1}{l}{X-n106-k14} &       & 2.0   & 26375.9 & \textbf{26362} &       & 4.0   & 26381.8 & 26378 &       & 17.6  & \textbf{26374.3} & \textbf{26362} &  \\
3     & \multicolumn{1}{l}{X-n110-k13} &       & 0.2   & \textbf{14971.0} & \textbf{14971} &       & 1.6   & \textbf{14971.0} & \textbf{14971} &       & 3.9   & \textbf{14971.0} & \textbf{14971} &  \\
4     & \multicolumn{1}{l}{X-n115-k10} &       & 0.2   & \textbf{12747.0} & \textbf{12747} &       & 1.8   & \textbf{12747.0} & \textbf{12747} &       & 6.4   & \textbf{12747.0} & \textbf{12747} &  \\
5     & \multicolumn{1}{l}{X-n120-k6} &       & 1.7   & 13337.6 & \textbf{13332} &       & 2.3   & \textbf{13332.0} & \textbf{13332} &       & 38.3  & \textbf{13332.0} & \textbf{13332} &  \\
6     & \multicolumn{1}{l}{X-n125-k30} &       & 1.4   & 55673.8 & \textbf{55539} &       & 2.7   & 55542.1 & \textbf{55539} &       & 6.1   & \textbf{55540.0} & \textbf{55539} &  \\
7     & \multicolumn{1}{l}{X-n129-k18} &       & 1.9   & 28998.0 & 28948 &       & 2.7   & 28948.5 & \textbf{28940} &       & 8.4   & \textbf{28940.0} & \textbf{28940} &  \\
8     & \multicolumn{1}{l}{X-n134-k13} &       & 2.1   & 10947.4 & \textbf{10916} &       & 3.3   & 10934.9 & \textbf{10916} &       & 20.2  & \textbf{10916.0} & \textbf{10916} &  \\
9     & \multicolumn{1}{l}{X-n139-k10} &       & 1.6   & 13603.1 & \textbf{13590} &       & 2.3   & \textbf{13590.0} & \textbf{13590} &       & 8.9   & \textbf{13590.0} & \textbf{13590} &  \\
\midrule
10    & \multicolumn{1}{l}{X-n143-k7} &       & 1.6   & 15745.2 & 15726 &       & 3.1   & 15700.2 & \textbf{15700} &       & 33.1  & \textbf{15700.0} & \textbf{15700} &  \\
11    & \multicolumn{1}{l}{X-n148-k46} &       & 0.8   & 43452.1 & \textbf{43448} &       & 3.2   & \textbf{43448.0} & \textbf{43448} &       & 5.1   & \textbf{43448.0} & \textbf{43448} &  \\
12    & \multicolumn{1}{l}{X-n153-k22} &       & 0.5   & 21400.0 & 21340 &       & 5.5   & 21226.3 & \textbf{21220} &       & 15.2  & \textbf{21225.6} & 21225 &  \\
13    & \multicolumn{1}{l}{X-n157-k13} &       & 0.8   & \textbf{16876.0} & \textbf{16876} &       & 3.2   & \textbf{16876.0} & \textbf{16876} &       & 28.4  & \textbf{16876.0} & \textbf{16876} &  \\
14    & \multicolumn{1}{l}{X-n162-k11} &       & 0.5   & 14160.1 & \textbf{14138} &       & 3.3   & 14141.3 & \textbf{14138} &       & 12.2  & \textbf{14138.0} & \textbf{14138} &  \\
15    & \multicolumn{1}{l}{X-n167-k10} &       & 0.9   & 20608.7 & 20562 &       & 3.7   & 20563.2 & \textbf{20557} &       & 36.1  & \textbf{20557.0} & \textbf{20557} &  \\
16    & \multicolumn{1}{l}{X-n172-k51} &       & 0.6   & 45616.1 & \textbf{45607} &       & 3.8   & \textbf{45607.0} & \textbf{45607} &       & 5.7   & \textbf{45607.0} & \textbf{45607} &  \\
17    & \multicolumn{1}{l}{X-n176-k26} &       & 1.1   & 48249.8 & 48140 &       & 7.6   & 47957.2 & \textbf{47812} &       & 16.2  & \textbf{47830.7} & \textbf{47812} &  \\
18    & \multicolumn{1}{l}{X-n181-k23} &       & 1.6   & 25571.5 & \textbf{25569} &       & 6.3   & 25591.1 & \textbf{25569} &       & 13.9  & \textbf{25569.4} & \textbf{25569} &  \\
19    & \multicolumn{1}{l}{X-n186-k15} &       & 1.7   & 24186.0 & \textbf{24145} &       & 5.9   & 24147.2 & \textbf{24145} &       & 20.2  & \textbf{24145.0} & \textbf{24145} &  \\
\midrule
20    & \multicolumn{1}{l}{X-n190-k8} &       & 2.1   & 17143.1 & 17085 &       & 12.1  & 16987.9 & \textbf{16980} &       & 161.9 & \textbf{16985.3} & \textbf{16980} &  \\
21    & \multicolumn{1}{l}{X-n195-k51} &       & 0.9   & \textbf{44234.3} & \textbf{44225} &       & 6.1   & 44244.1 & \textbf{44225} &       & 9.3   & 44283.8 & \textbf{44225} &  \\
22    & \multicolumn{1}{l}{X-n200-k36} &       & 7.5   & 58697.2 & 58626 &       & 8.0   & 58626.4 & \textbf{58578} &       & 12.0  & \textbf{58615.1} & \textbf{58578} &  \\
23    & \multicolumn{1}{l}{X-n204-k19} &       & 1.1   & 19625.2 & 19570 &       & 5.4   & 19571.5 & \textbf{19565} &       & 15.7  & \textbf{19567.0} & \textbf{19565} &  \\
24    & \multicolumn{1}{l}{X-n209-k16} &       & 3.8   & 30765.4 & 30667 &       & 8.6   & 30680.4 & \textbf{30656} &       & 35.7  & \textbf{30671.3} & \textbf{30656} &  \\
25    & \multicolumn{1}{l}{X-n214-k11} &       & 2.3   & 11126.9 & 10985 &       & 10.2  & 10877.4 & \textbf{10856} &       & 52.3  & \textbf{10872.1} & \textbf{10856} &  \\
26    & \multicolumn{1}{l}{X-n219-k73} &       & 0.9   & \textbf{117595.0} & \textbf{117595} &       & 7.7   & 117604.9 & \textbf{117595} &       & 18.2  & 117600.5 & \textbf{117595} &  \\
27    & \multicolumn{1}{l}{X-n223-k34} &       & 8.5   & 40533.5 & 40471 &       & 8.3   & 40499.0 & \textbf{40437} &       & 18.6  & \textbf{40478.4} & \textbf{40437} &  \\
28    & \multicolumn{1}{l}{X-n228-k23} &       & 2.4   & 25795.8 & 25743 &       & 9.8   & 25779.3 & \textbf{25742} &       & 29.0  & \textbf{25768.0} & 25743 &  \\
29    & \multicolumn{1}{l}{X-n233-k16} &       & 3.0   & 19336.7 & 19266 &       & 6.8   & 19288.4 & \textbf{19230} &       & 34.0  & \textbf{19276.5} & \textbf{19230} &  \\
\midrule
30    & \multicolumn{1}{l}{X-n237-k14} &       & 3.5   & 27078.8 & \textbf{27042} &       & 8.9   & 27067.3 & \textbf{27042} &       & 32.2  & \textbf{27048.8} & \textbf{27042} &  \\
31    & \multicolumn{1}{l}{X-n242-k48} &       & 17.8  & \textbf{82874.2} & 82774 &       & 12.4  & 82948.7 & 82804 &       & 18.1  & 82920.9 & \textbf{82751} & \\
32    & \multicolumn{1}{l}{X-n247-k47} &       & 2.1   & 37507.2 & 37289 &       & 20.4  & \textbf{37284.4} & \textbf{37274} &       & 27.7  & 37388.9 & \textbf{37274} &  \\
33    & \multicolumn{1}{l}{X-n251-k28} &       & 10.8  & 38840.0 & 38727 &       & 11.7  & 38796.4 & 38699 &       & 20.2  & \textbf{38778.7} & \textbf{38684} & \\
34    & \multicolumn{1}{l}{X-n256-k16} &       & 2.0   & 18883.9 & 18880 &       & 6.5   & 18880.0 & 18880 &       & 23.0  & \textbf{18867.7} & \textbf{18839} & $\circledast$ \\
35    & \multicolumn{1}{l}{X-n261-k13} &       & 6.7   & 26869.0 & 26706 &       & 12.7  & 26629.6 & \textbf{26558} &       & 48.6  & \textbf{26618.1} & \textbf{26558} &  \\
36    & \multicolumn{1}{l}{X-n266-k58} &       & 10.0  & \textbf{75563.3} & \textbf{75478} &       & 21.4  & 75759.3 & 75517 &       & 29.9  & 75710.7 & \textbf{75478} &  \\
37    & \multicolumn{1}{l}{X-n270-k35} &       & 9.1   & 35363.4 & 35324 &       & 11.3  & 35367.2 & \textbf{35303} &       & 18.9  & \textbf{35314.6} & \textbf{35303} &  \\
38    & \multicolumn{1}{l}{X-n275-k28} &       & 3.6   & 21256.0 & \textbf{21245} &       & 12.0  & 21280.6 & \textbf{21245} &       & 22.7  & \textbf{21255.0} & \textbf{21245} &  \\
39    & \multicolumn{1}{l}{X-n280-k17} &       & 9.6   & 33769.4 & 33624 &       & 19.1  & 33605.8 & 33505 &       & 136.2 & \textbf{33587.9} & \textbf{33503} & \\
\midrule
40    & \multicolumn{1}{l}{X-n284-k15} &       & 8.6   & 20448.5 & 20295 &       & 19.9  & 20286.4 & \textbf{20227} &       & 97.7  & \textbf{20282.1} & 20228 &  \\
41    & \multicolumn{1}{l}{X-n289-k60} &       & 16.1  & 95450.6 & 95315 &       & 21.3  & 95469.5 & 95244 &       & 41.7  & \textbf{95447.2} & \textbf{95211} &  \\
42    & \multicolumn{1}{l}{X-n294-k50} &       & 12.4  & \textbf{47254.7} & 47190 &       & 14.7  & 47259.0 & 47171 &       & 27.0  & 47272.7 & \textbf{47161} & $\circledast$ \\
43    & \multicolumn{1}{l}{X-n298-k31} &       & 6.9   & 34356.0 & 34239 &       & 10.9  & 34292.1 & \textbf{34231} &       & 20.7  & \textbf{34276.3} & \textbf{34231} &  \\
44    & \multicolumn{1}{l}{X-n303-k21} &       & 14.2  & 21895.8 & 21812 &       & 17.3  & 21850.9 & 21748 &       & 48.4  & \textbf{21811.2} & \textbf{21744} &  \\
45    & \multicolumn{1}{l}{X-n308-k13} &       & 9.5   & 26101.1 & 25901 &       & 15.3  & \textbf{25895.4} & \textbf{25859} &       & 112.8 & 25897.3 & 25861 &  \\
46    & \multicolumn{1}{l}{X-n313-k71} &       & 17.5  & 94297.3 & 94192 &       & 22.4  & \textbf{94265.2} & 94093 &       & 30.6  & 94280.4 & \textbf{94045} & \\
47    & \multicolumn{1}{l}{X-n317-k53} &       & 8.6   & \textbf{78356.0} & \textbf{78355} &       & 22.4  & 78387.8 & \textbf{78355} &       & 50.3  & 78385.3 & \textbf{78355} &  \\
48    & \multicolumn{1}{l}{X-n322-k28} &       & 14.7  & 29991.3 & 29877 &       & 15.2  & 29956.1 & 29870 &       & 27.7  & \textbf{29892.5} & \textbf{29834} & $\circledast$ \\
49    & \multicolumn{1}{l}{X-n327-k20} &       & 19.1  & 27812.4 & 27599 &       & 18.2  & 27628.2 & 27564 &       & 68.7  & \textbf{27590.8} & \textbf{27532} & $\circledast$ \\
50    & \multicolumn{1}{l}{X-n331-k15} &       & 15.7  & 31235.5 & 31105 &       & 24.4  & 31159.6 & \textbf{31103} &       & 102.1 & \textbf{31126.7} & \textbf{31103} &  \\
\midrule
\multicolumn{2}{r}{Average gap:} &       &       & 0.37\% & 0.13\% &       &       & 0.14\% & 0.02\% &       &       & 0.10\% & 0.00\% &  \\
\bottomrule
\end{tabular}%

    \end{table}

    \begin{table}[htbp]\scriptsize
      \centering
      \caption{Results for the instances with more than 331 customers from \cite{Uchoa2017}}
      \label{cvrp:tab:large}
      \def\arraystretch{1.05}\tabcolsep=3.0pt
\begin{tabular}{rrrrrrrrrrrrrrrr}
\toprule
\multirow{2}[4]{*}{\#} & \multicolumn{1}{l}{\multirow{2}[4]{*}{Instance}} &       & \multicolumn{3}{c}{ILS} &       & \multicolumn{3}{c}{UHGS} &       & \multicolumn{4}{c}{UHGS-BS} \\
\cmidrule{4-6}\cmidrule{8-10}\cmidrule{12-15}      &       &       & Time  & Average & Best  &       & Time  & Average & Best  &       & Time  & Average & Best  &   \\
\midrule
51    & \multicolumn{1}{l}{X-n336-k84} &       & 21.4  & 139461.0 & \textbf{139197} &       & 38.0  & 139534.9 & 139210 &       & 66.0  & \textbf{139460.1} & 139303 &  \\
52    & \multicolumn{1}{l}{X-n344-k43} &       & 22.6  & 42284.0 & 42146 &       & 21.7  & 42208.8 & 42099 &       & 39.7  & \textbf{42156.1} & \textbf{42056} & $\circledast$ \\
53    & \multicolumn{1}{l}{X-n351-k40} &       & 25.2  & 26150.3 & 26021 &       & 33.7  & 26014.0 & 25946 &       & 51.5  & \textbf{25981.8} & \textbf{25938} &  \\
54    & \multicolumn{1}{l}{X-n359-k29} &       & 48.9  & 52076.5 & 51706 &       & 34.9  & 51721.7 & \textbf{51509} &       & 112.0 & \textbf{51640.7} & 51555 &  \\
55    & \multicolumn{1}{l}{X-n367-k17} &       & 13.1  & 23003.2 & 22902 &       & 22.0  & \textbf{22838.4} & \textbf{22814} &       & 117.3 & 22876.2 & \textbf{22814} &  \\
56    & \multicolumn{1}{l}{X-n376-k94} &       & 7.1   & \textbf{147713.0} & \textbf{147713} &       & 28.3  & 147750.2 & 147717 &       & 70.3  & 147740.5 & 147714 &  \\
57    & \multicolumn{1}{l}{X-n384-k52} &       & 34.5  & 66372.5 & 66116 &       & 40.2  & 66270.2 & 66081 &       & 56.8  & \textbf{66170.3} & \textbf{65997} &  \\
58    & \multicolumn{1}{l}{X-n393-k38} &       & 20.8  & 38457.4 & 38298 &       & 28.7  & 38374.9 & 38269 &       & 49.3  & \textbf{38309.3} & \textbf{38260} & $\circledast$ \\
59    & \multicolumn{1}{l}{X-n401-k29} &       & 60.4  & 66715.1 & 66453 &       & 49.5  & 66365.4 & 66243 &       & 110.2 & \textbf{66359.0} & \textbf{66212} & \\
\midrule
60    & \multicolumn{1}{l}{X-n411-k19} &       & 23.8  & 19954.9 & 19792 &       & 34.7  & 19743.8 & \textbf{19718} &       & 126.0 & \textbf{19736.7} & 19721 &  \\
61    & \multicolumn{1}{l}{X-n420-k130} &       & 22.2  & \textbf{107838.0} & \textbf{107798} &       & 53.2  & 107924.1 & \textbf{107798} &       & 87.7  & 107913.7 & \textbf{107798} &  \\
62    & \multicolumn{1}{l}{X-n429-k61} &       & 38.2  & 65746.6 & 65563 &       & 41.5  & \textbf{65648.5} & 65501 &       & 65.6  & 65661.6 & \textbf{65470} & \\
63    & \multicolumn{1}{l}{X-n439-k37} &       & 39.6  & 36441.6 & \textbf{36395} &       & 34.6  & 36451.1 & \textbf{36395} &       & 57.1  & \textbf{36410.1} & \textbf{36395} &  \\
64    & \multicolumn{1}{l}{X-n449-k29} &       & 59.9  & 56204.9 & 55761 &       & 64.9  & 55553.1 & 55378 &       & 132.6 & \textbf{55432.7} & \textbf{55330} &  \\
65    & \multicolumn{1}{l}{X-n459-k26} &       & 60.6  & 24462.4 & 24209 &       & 42.8  & 24272.6 & 24181 &       & 92.9  & \textbf{24226.0} & \textbf{24145} & $\circledast$ \\
66    & \multicolumn{1}{l}{X-n469-k138} &       & 36.3  & \textbf{222182.0} & \textbf{221909} &       & 86.7  & 222617.1 & 222070 &       & 142.3 & 222427.5 & 222235 &  \\
67    & \multicolumn{1}{l}{X-n480-k70} &       & 50.4  & 89871.2 & 89694 &       & 67.0  & 89760.1 & 89535 &       & 73.1  & \textbf{89744.7} & \textbf{89513} &  \\
68    & \multicolumn{1}{l}{X-n491-k59} &       & 52.2  & 67226.7 & 66965 &       & 71.9  & 66898.0 & 66633 &       & 81.9  & \textbf{66794.1} & \textbf{66607} &  \\
69    & \multicolumn{1}{l}{X-n502-k39} &       & 80.8  & 69346.8 & 69284 &       & 63.6  & 69328.8 & 69253 &       & 177.7 & \textbf{69277.1} & \textbf{69247} &  \\
\midrule
70    & \multicolumn{1}{l}{X-n513-k21} &       & 35.0  & 24434.0 & 24332 &       & 33.1  & 24296.6 & \textbf{24201} &       & 99.4  & \textbf{24256.2} & \textbf{24201} &  \\
71    & \multicolumn{1}{l}{X-n524-k137} &       & 27.3  & 155005.0 & \textbf{154709} &       & 80.7  & \textbf{154979.5} & 154774 &       & 207.3 & 155038.1 & 154787 &  \\
72    & \multicolumn{1}{l}{X-n536-k96} &       & 62.1  & 95700.7 & 95524 &       & 107.5 & \textbf{95330.6} & 95122 &       & 144.5 & 95335.4 & \textbf{95112} &  \\
73    & \multicolumn{1}{l}{X-n548-k50} &       & 64.0  & \textbf{86874.1} & \textbf{86710} &       & 84.2  & 86998.5 & 86822 &       & 136.6 & 86881.0 & 86778 &  \\
74    & \multicolumn{1}{l}{X-n561-k42} &       & 68.9  & 43131.3 & 42952 &       & 60.6  & 42866.4 & 42756 &       & 77.2  & \textbf{42860.0} & \textbf{42733} &  \\
75    & \multicolumn{1}{l}{X-n573-k30} &       & 112.0 & 51173.0 & 51092 &       & 188.2 & 50915.1 & \textbf{50780} &       & 782.4 & \textbf{50876.9} & 50801 &  \\
76    & \multicolumn{1}{l}{X-n586-k159} &       & 78.5  & 190919.0 & 190612 &       & 175.3 & 190838.0 & 190543 &       & 234.3 & \textbf{190752.4} & \textbf{190442} &  \\
77    & \multicolumn{1}{l}{X-n599-k92} &       & 73.0  & 109384.0 & 109056 &       & 125.9 & 109064.2 & 108813 &       & 166.9 & \textbf{108993.3} & \textbf{108576} &  \\
78    & \multicolumn{1}{l}{X-n613-k62} &       & 74.8  & 60444.2 & 60229 &       & 117.3 & 59960.0 & 59778 &       & 103.6 & \textbf{59859.7} & \textbf{59654} &  \\
79    & \multicolumn{1}{l}{X-n627-k43} &       & 162.7 & 62905.6 & 62783 &       & 239.7 & 62524.1 & 62366 &       & 543.1 & \textbf{62442.9} & \textbf{62254} &  \\
\midrule
80    & \multicolumn{1}{l}{X-n641-k35} &       & 140.4 & 64606.1 & 64462 &       & 158.8 & 64192.0 & \textbf{63839} &       & 304.4 & \textbf{64105.6} & 63859 &  \\
81    & \multicolumn{1}{l}{X-n655-k131} &       & 47.2  & \textbf{106782.0} & \textbf{106780} &       & 150.5 & 106899.1 & 106829 &       & 253.2 & 106855.6 & 106804 &  \\
82    & \multicolumn{1}{l}{X-n670-k126} &       & 61.2  & 147676.0 & 147045 &       & 264.1 & \textbf{147222.7} & \textbf{146705} &       & 267.7 & 147663.9 & 147163 &  \\
83    & \multicolumn{1}{l}{X-n685-k75} &       & 73.9  & 68988.2 & 68646 &       & 156.7 & 68654.1 & \textbf{68425} &       & 177.0 & \textbf{68596.0} & 68496 &  \\
84    & \multicolumn{1}{l}{X-n701-k44} &       & 210.1 & 83042.2 & 82888 &       & 253.2 & 82487.4 & 82293 &       & 368.0 & \textbf{82409.2} & \textbf{82174} &  \\
85    & \multicolumn{1}{l}{X-n716-k35} &       & 225.8 & 44171.6 & 44021 &       & 264.3 & 43641.4 & 43525 &       & 437.2 & \textbf{43599.9} & \textbf{43498} &  \\
86    & \multicolumn{1}{l}{X-n733-k159} &       & 111.6 & 137045.0 & 136832 &       & 244.5 & \textbf{136587.6} & \textbf{136366} &       & 334.2 & 136607.4 & 136424 &  \\
87    & \multicolumn{1}{l}{X-n749-k98} &       & 127.2 & 78275.9 & 77952 &       & 313.9 & 77864.9 & 77715 &       & 308.3 & \textbf{77862.8} & \textbf{77605} &  \\
88    & \multicolumn{1}{l}{X-n766-k71} &       & 242.1 & 115738.0 & 115443 &       & 383.0 & 115147.9 & \textbf{114683} &       & 330.5 & \textbf{115115.9} & 114812 &  \\
89    & \multicolumn{1}{l}{X-n783-k48} &       & 235.5 & 73722.9 & 73447 &       & 269.7 & 73009.6 & 72781 &       & 351.2 & \textbf{72892.4} & \textbf{72738} &  \\
\midrule
90    & \multicolumn{1}{l}{X-n801-k40} &       & 432.6 & 74005.7 & 73830 &       & 289.2 & 73731.0 & 73587 &       & 424.0 & \textbf{73651.6} & \textbf{73466} &  \\
91    & \multicolumn{1}{l}{X-n819-k171} &       & 148.9 & 159425.0 & 159164 &       & 374.3 & 158899.3 & 158611 &       & 675.6 & \textbf{158849.0} & \textbf{158592} &  \\
92    & \multicolumn{1}{l}{X-n837-k142} &       & 173.2 & 195027.0 & 194804 &       & 463.4 & \textbf{194476.5} & \textbf{194266} &       & 634.9 & 194504.0 & 194356 &  \\
93    & \multicolumn{1}{l}{X-n856-k95} &       & 153.7 & 89277.6 & 89060 &       & 288.4 & 89238.7 & 89118 &       & 314.6 & \textbf{89220.0} & \textbf{89020} &  \\
94    & \multicolumn{1}{l}{X-n876-k59} &       & 409.3 & 100417.0 & 100177 &       & 495.4 & 99884.1 & 99715 &       & 543.1 & \textbf{99780.3} & \textbf{99610} &  \\
95    & \multicolumn{1}{l}{X-n895-k37} &       & 410.2 & 54958.5 & 54713 &       & 321.9 & 54439.8 & \textbf{54172} &       & 500.2 & \textbf{54407.4} & 54254 &  \\
96    & \multicolumn{1}{l}{X-n916-k207} &       & 226.1 & 330948.0 & 330639 &       & 560.8 & 330198.3 & \textbf{329836} &       & 1082.5 & \textbf{330153.2} & 329866 &  \\
97    & \multicolumn{1}{l}{X-n936-k151} &       & 202.5 & 134530.0 & 133592 &       & 531.5 & \textbf{133512.9} & \textbf{133140} &       & 1022.2 & 133729.3 & 133376 &  \\
98    & \multicolumn{1}{l}{X-n957-k87} &       & 311.2 & 85936.6 & 85697 &       & 432.9 & 85822.6 & 85672 &       & 307.9 & \textbf{85681.5} & \textbf{85555} &  \\
99    & \multicolumn{1}{l}{X-n979-k58} &       & 687.2 & 120253.0 & 119994 &       & 554.0 & \textbf{119502.1} & 119194 &       & 928.4 & 119527.7 & \textbf{119188} &  \\
100   & \multicolumn{1}{l}{X-n1001-k43} &       & 792.8 & 73985.4 & 73776 &       & 549.0 & 72956.0 & 72742 &       & 952.8 & \textbf{72903.3} & \textbf{72629} &  \\
\midrule
\multicolumn{2}{r}{Average gap:} &       &       & 0.74\% & 0.42\% &       &       & 0.30\% & 0.06\% &       &       & 0.24\% & 0.03\% &  \\
\bottomrule
\end{tabular}%

    \end{table}

These tables highlight how the local search applied on $\SB{2}$ resulted in several improvements over the best solutions obtained by the state-of-the-art methods considered, with an average gap of 0.10\% and 0.24\% on the medium and large instances, respectively, in comparison with 0.14\% and 0.30\% for a classical search in $\cS$. This improvement in solution quality comes at the price of an overall twofold increase of CPU time.

Another notable observation of these experiments is that the search in $\SB{2}$ finds solutions which are structurally different. Indeed, we obtained some new best known solutions for unexpectedly small instances, with 256, 294, 322 and 327 customers (marked with a $\circledast$). This is not a coincidence, given that the BKS listed at \url{http://vrp.atd-lab.inf.puc-rio.br/} already originate from various previous articles and methods, over a large cumulated amount of test runs and parameter settings. For the particular case of the instance X-n256-k16, one interesting characteristic was observed: only 16 vehicles are used, with a total capacity usage of 99.6\%.

\FloatBarrier
\section{Conclusions and future work}
\label{cvrp:sec:conclusion}
 
In this article, we investigated decision-set decompositions for the classic CVRP. Our experiments show that decomposing the problem into \textsc{Assignment} and \textsc{Sequencing} decisions, and conducting local search in the \textsc{Assignment} space ($\SA$) generates consistently better results than heuristically searching in the complete search space $\cS$. When doing so, each solution is systematically \emph{decoded} by finding optimal routes (\textsc{Sequencing} decisions) with the \textsc{Concorde} TSP solver. However, the extra CPU dedicated to solution decoding is prohibitively high to employ this technique within state-of-the-art metaheuristics. To circumvent this issue, the $\BS$ neighborhood was employed to define an intermediate search space ($\SB{k}$) which can be more efficiently searched while retaining a high solution quality.

Different techniques were proposed and evaluated for efficiently searching in space $\SB{k}$: \emph{neighborhood reduction}, \emph{dynamic move filters}, \emph{concatenation techniques} and \emph{efficient memory structures}. Moreover, tunneling techniques were employed to reshape $\SB{k}$ into a search space more similar to $\SA$ as the search progresses.
The combination of these techniques within the UHGS solver resulted in an significant improvement of solution accuracy. Multiple instances from the literature had their best known solution improved, and new state-of-the-art results were obtained for the CVRP.

This improvement of search space, however, still results in some extra computational effort.
Therefore, many research perspectives are open about how to fully exploit a larger search space such as $\SB{k}$ without any time compromise. One possibility could be to employ the search in $\SB{k}$ only in exceptional circumstances, for a selected subset of promising solutions (e.g., each new best incumbent solutions during the search). Other open research possibilities relate to the search space $\SA$. Indeed, even with efficient memory structures and neighborhood reduction techniques, using \text{Concorde} for each solution evaluation remains impracticable. This effort could be mitigated if good and fast lower bounds were proposed for the cost of the routes, therefore permitting to filter a large proportion of moves as in \cite{Vidal2015a}. \textsc{Concorde} is also not optimized to handle millions of small cardinality routes issued from a local search, such that a dedicated TSP solution procedure exploiting information from the current incumbent solution represents another promising research line. Finally, the proposed approach can be naturally evaluated for other variants of the CVRP, in the presence of different types of attributes that have an impact on the \textsc{Assignment} and \textsc{Sequencing} decision classes. These are all promising perspectives for future research at the crossroads of dynamic programming, integer programming, and metaheuristics.

\section*{Acknowledgements}
\label{sec:acknowledgements}

This research has been funded by the Belgian Science Policy Office (BELSPO) in the Interuniversity Attraction Pole COMEX (http://comex.ulb.ac.be), the National Council for Scientific and Technological Development (CNPQ -- grant number 308498/2015-1) and FAPERJ in Brazil (grant number E-26/203.310/2016). The computational resources and services used in this work were provided by the VSC (Flemish Supercomputer Center), funded by the Hercules Foundation and the Flemish Government -- department EWI. This support is gratefully acknowledged. In addition, we would like to thank Luke Connolly for providing editorial consultation.


\end{document}